\documentclass[12pt]{amsart}
\usepackage{amssymb}
\usepackage{graphics}
\input xy
\xyoption{all}

\newtheorem{theorem}{Theorem}[section]
\newtheorem{corollary}[theorem]{Corollary}
\newtheorem{lemma}[theorem]{Lemma}
\newtheorem{proposition}[theorem]{Proposition}
\newtheorem{question}[theorem]{Question}
\newtheorem{problem}[theorem]{Problem}

\theoremstyle{definition}
\newtheorem{definition}[theorem]{Definition}

\theoremstyle{remark}
\newtheorem{remark}[theorem]{Remark}
\newtheorem{rem}[theorem]{Remark}
\numberwithin{equation}{section}

\newcommand{\GL}{\mathrm{GL}}
\newcommand{\Mod}{\mathrm{Mod}}

\newcommand{\Homeo}{\mathrm{Homeo}}
\newcommand{\Core}{\mathrm{Core}}
\newcommand{\Shadow}{\mathrm{Shadow}}
\newcommand{\Conv}{\mathrm{Conv}}
\newcommand{\Fix}{\mathrm{Fix}}
\newcommand{\CAT}{\mathrm{CAT}}
\newcommand{\Aff}{\mathrm{Aff}}
\newcommand{\PMF}{\mathcal{PMF}}
\newcommand{\SL}{\mathrm{SL}}
\newcommand{\PSL}{\mathrm{PSL}}
\newcommand{\PGL}{\mathrm{PGL}}
\newcommand{\Ad}{\mathrm{Ad}}
\newcommand{\Aut}{\mathrm{Aut}}
\newcommand{\Inn}{\mathrm{Inn}}
\newcommand{\Out}{\mathrm{Out}}

\newcommand{\arrow}{\rightarrow}
\newcommand{\trivgp}{\langle e \rangle}
\newcommand{\defeq}{\stackrel{\mathrm{def}}{=}}
\newcommand{\R}{\mathbb R}
\newcommand{\C}{\mathbb C}
\newcommand{\Q}{\mathbb Q}
\newcommand{\N}{\mathbb N}
\newcommand{\Z}{\mathbb Z}
\newcommand{\F}{\mathbb F}

\newcommand{\PP}{\mathbb P}
\newcommand{\G}{\mathbb G}
\renewcommand{\H}{\mathbb H}

\newcommand{\gC}{\Gamma}\newcommand{\gS}{\Sigma}
\newcommand{\gD}{\Delta}

\newcommand{\gd}{\delta}
\newcommand{\gep}{\epsilon}
\newcommand{\gb}{\beta}
\newcommand{\ga}{\alpha}

\newcommand{\ti}{\tilde}

\begin{document}
\title{Countable primitive groups.}%
\author{Tsachik Gelander}%
\address{Yale University.}%
\email{tsachik.gelander@yale.edu}%

\author{Yair Glasner}
\address{Institute for advanced study.}
\email{yair@math.uic.edu}
\date{\today}

\thanks{T.G. partially supported by NSF grant
DMS-0404557, Y.G. was supported by NSF grant DMS-0111298}

\begin{abstract}
We give a complete characterization of countable primitive groups in several settings including
linear groups, subgroups of mapping class groups, groups acting minimally on trees and convergence
groups. The latter category includes as a special case Kleinian groups as well as subgroups of
word hyperbolic groups. As an application we calculate the Frattini subgroup in many of these
settings, often generalizing results that were only known for finitely generated groups. In
particular, we answer a question of G. Higman and B.H. Neumann on the Frattini group of an
amalgamated product.
\end{abstract}

\maketitle
%\tableofcontents

% ----------------------------------------------------------------
%\subjclass[2000]{Primary 20B15; Secondary 20B07,20E28}%
%\keywords{}%

%\date{\today}

\section{Introduction} \label{sec:introduction}

In the study of groups one often tries to understand their properties via actions on various
geometric, combinatorial or algebraic structures. The most fundamental mathematical structure that
comes to mind is a set. When groups were first introduced by Galois, they were considered as
permutation groups (of roots of polynomials). The abstract definition of a group, without any
given realization as a symmetry group appeared only later. Today permutation representations are
indispensable in the study of finite groups. There is a beautiful theory of infinite permutation
groups, which has elaborate connections with logic and model theory. But this theory does not
often make it easier to understand a given infinite group. In this paper we establish connections
between permutation representation theory to other representation theories.

\subsection{Terminology}

\begin{definition}
An action of a group $\Gamma$ on a set $X$ is {\it primitive} if $|X|>1$ and there are no
$\Gamma$-invariant equivalence relations on $X$ apart from the two trivial ones\footnote{The
trivial equivalence relations are those with a unique equivalence class, or with singletons as
equivalence classes. When $|X|=2$, one should also require that the action is not trivial.}. An
action is called {\it quasiprimitive} if every normal subgroup acts either trivially or
transitively. A group is {\it primitive} or {\it quasiprimitive} if it admits a faithful primitive
or quasiprimitive action on a set.
\end{definition}

Primitive actions are the irreducible building blocks in the theory of permutation
representations: any imprimitive action can be embedded into a wreath product of two simpler ones
-- the action on the set of equivalence classes, and the action on a given equivalence class by
its setwise stabilizer.

Any action that is transitive on pairs is also primitive. A primitive action, in turn, is
quasiprimitive because orbits of normal subgroups define invariant equivalence relations. All
these actions are in particular transitive, and the above properties can be characterized in terms
of the stabilizer of a point. The action $\Gamma \circlearrowleft \Gamma/\Delta$ is primitive if
and only if $\Delta < \Gamma$ is maximal. The kernel of this action or the {\it core of $\Delta$
in $\Gamma$} is defined by $\Core_{\Gamma}(\Delta) = \cap_{\gamma \in \Gamma}(\Delta^{\gamma}) =
\trivgp$, where $\Delta ^{\gamma} = \gamma \Delta \gamma^{-1}$. In particular the action is
faithful if and only if $\Core_{\Gamma}(\Delta) = \trivgp$. Finally the action on $\Gamma
\circlearrowleft \Gamma/\Delta$ is quasiprimitive and faithful if and only if $\Delta$ is a proper
prodense subgroup of $\Gamma$ in the following sense.
\begin{definition} \label{def:prodense}
A subgroup $\Delta \leq \Gamma$ is called prodense if $\Delta N = \Gamma$ for every non-trivial
normal subgroup $\trivgp \ne N \lhd \Gamma.$
\end{definition}

\noindent Assume that the intersection of any two non-trivial normal subgroups of $\Gamma$ is
still non-trivial. Then the collection of cosets of non-trivial normal subgroups forms a basis for
an invariant topology on the group, which we call {\it the normal topology}. In that case, a
subgroup is prodense if and only if it is dense in this topology. This is the origin of the name
prodense. It is interesting to note that for residually finite quasiprimitive groups the normal
topology is always well defined (c.f. Corollary \ref{cor:RF_non_banal}).

\subsection{Goals.}
Our goal is to establish connections between the structure of a group and its permutation
representation theory. We address the following

\begin{problem} \label{Q:main}
Characterize primitive groups -- the groups that admit a faithful primitive action on a set.
\end{problem}

\noindent While this question in complete generality seems hopeless, we develop tools that enable
us to answer this question in certain {\it geometric settings}: linear groups, subgroups of
mapping class groups, groups acting minimally on trees, and convergence groups. In each one of
these geometric setting we obtain a complete characterization of the countable primitive groups.
The criteria for primitivity are usually explicit and easy to check, and hence produce various
examples and counterexamples.

\subsection{About minimal normal subgroups}
\label{sec:nt_commuting} Minimal normal subgroups play a central role in the study of finite
primitive actions. The novelty of this paper lies in abandoning this approach and appealing to the
normal topology instead. When the normal topology is not well defined, or more generally when the
group is banal in the sense of the following definition, our methods fail.
\begin{definition} \label{def:banal}
We say that a quasiprimitive group $\Gamma$ is {\it banal} if there exist non-trivial normal
subgroups $\trivgp \ne M, N \lhd \Gamma$ that commute elementwise $[M,N] = \trivgp$.
\end{definition}

Luckily, banal groups always contain minimal normal subgroups, so classical methods can be used to
understand their primitive action. This analysis yields exactly two types of banal groups which
are defined as follows:

\begin{definition}\label{def:affine}
Let $M$ be a vector space over a prime field, and let $\Delta \leq \GL (M)$ be such that there are
no non-trivial $\Delta$ invariant subgroups of (the additive group of) $M$. The group $\Gamma =
\Delta \ltimes M$ admits a natural {\it affine} action on $M$, where $M$ acts on itself by (left)
translation, and $\Delta$ acts by conjugation. We say that the permutation group $\Gamma$ is {\it
primitive of affine type}.
\end{definition}

\begin{definition}\label{def:diagonal}
Let $M$ be a nonabelian characteristically simple group, and $\Delta\leq\Aut(M)$ a subgroup
containing $M \cong \Inn(M)$ such that there are no non-trivial $\Delta$-invariant subgroups of
$M$. The group $\Gamma = \Delta \ltimes M$ admits a natural affine action on $M$ as above. We say
that the permutation group $\Gamma$ is {\it primitive of diagonal type}\footnote{Note that a group
$\Gamma$ of diagonal type as above, contains another normal subgroup isomorphic to $M$, namely
$\{i(m^{-1}) m \ | \ m \in M \}$ where $i: M \arrow \Inn(M) < \Delta$ is the natural injection.}.
\end{definition}

It is not difficult to verify that the affine action of these groups is indeed primitive.

As we are mostly interested in geometric methods that do not involve minimal normal subgroups, we
postpone the analysis of banal groups to the appendix, and assume throughout the paper that our
groups are not banal. A posteriori this approach is justified, because it turns out that the
non-existence of minimal normal subgroups is quite typical. In fact the only time that banal
groups actually appear in this paper is when linear groups that are not finitely generated are
treated.

\subsection{The Margulis--So\u{\i}fer theorem}
Let us mention the beautiful theorem of Margulis and So\u{\i}fer which completely characterizes
finitely generated linear groups that admit infinite primitive actions\footnote{In particular it
answered a question of Platonov about the existence of a maximal subgroup of infinite index in
$\SL_3(\Z )$.} and which was a main inspiration for our work.

\begin{theorem}\label{MS} (Margulis and So\u{\i}fer \cite{MS:Maximal})
A finitely generated linear group admits a primitive action on an infinite set if and only if it
is not virtually solvable.
\end{theorem}

Notwithstanding the elementary and simple formulation of Theorem \ref{MS}, the proof requires a
deep understanding of linear groups. However, the actions that Theorem \ref{MS} and its proof
provide are mostly non faithful. For example, one sees from the statement of \ref{MS} that the
property of admitting an infinite primitive action is stable (for such groups) under
commensurability. On the other hand Theorem \ref{thm:Linear} below allows one to construct
examples of primitive linear groups with finite index subgroups and supergroups which are not
primitive -- this illustrates the sensitivity of that stronger property.

\subsection{Statements of the main results}

Before listing our main results, let us note that in all the cases under consideration, we prove
that a group is primitive if and only if it is quasiprimitive.

\begin{definition} \label{def:lin_cond}
We say that a group $\Gamma$ satisfies the {\it linear conditions for primitivity} if it admits a
faithful linear representation over some algebraically closed field $\rho: \Gamma \to \GL_n(k)$
with Zariski closure $\G = \overline{\rho(\Gamma)}^{Z}$ satisfying the following conditions:
\begin{enumerate}
\item[{\rm (i)}] $\G^{\circ}$, the identity component, decomposes
as a direct product of simple factors with trivial centers.

\item[{\rm (ii)}] The action of $\Gamma$, by conjugation, on
$\G^{\circ}$ is faithful and permutes the simple factors of $\G^{\circ}$ transitively.
\end{enumerate}
\end{definition}

\noindent In particular any linear group with a simple Zariski closure (e.g. a lattice in a simple
center-free Lie group) satisfies the linear conditions for primitivity.

Our first theorem can be considered as a coarse generalization of O'Nan-Scott theorem (see
\cite{AS:Maximal_subgroups,DM:Permutation_Groups}) to the setting of countable linear groups.
\begin{theorem}\label{thm:Linear}
A countable non torsion linear group $\Gamma$ is primitive if and only if one of the following
mutually exclusive conditions hold.
\begin{itemize}
\item \label{itm:ZClosure} $\Gamma$ satisfies the linear conditions for primitivity as in
Definition \ref{def:lin_cond}.
\item \label{itm:affine} $\Gamma$  is primitive of affine type as in Definition \ref{def:affine}
\item \label{itm:diagonal} $\Gamma$ is primitive of diagonal type as in Definition
\ref{def:diagonal}.
\end{itemize}
In the affine and the diagonal cases the group $\Gamma$ is banal and it admits a unique
quasiprimitive action. For a finitely generated group $\Gamma$ only the first possibility
\ref{itm:ZClosure} can occur.
\end{theorem}

\begin{rem} In the 0 characteristic case, the theorem remains
valid without the assumption that the group is non-torsion. In
positive characteristic, we need this assumption for our proof. In
fact, in the proof of Theorem \ref{thm:Linear}, we actually
establish a stronger statement: the existence of a free prodense
subgroup which is contained in a maximal subgroup. This stronger
statement fails for torsion groups like $\PSL_2(\overline{\F_7})$,
where $\overline{\F_7}$ is the algebraic closure of $\F_7$. Note
however that $\PSL_2(\overline{\F_7})$ does not violate Theorem
\ref{thm:Linear} because it is primitive, and in fact even admits
a faithful 3-transitive action on the projective line
$\mathbb{P}\overline{F_7}$.
\end{rem}

As a corollary we can prove the following generalization of the Margulis-So\u{\i}fer theorem to
countable linear groups which are not necessarily finitely generated.

\begin{corollary}\label{GMS}
A countable linear non-torsion group which is not virtually solvable has a maximal subgroup of
infinite index.
\end{corollary}

\begin{rem} The condition that $\Gamma$ is not virtually solvable in Corollary \ref{GMS} is sufficient but no longer necessary.
For example the solvable group $\Aff(\Q)
\defeq \Q^{*} \ltimes \Q$ acts two transitively on the affine line
$\Q$, and in particular it is primitive (of affine type).
\end{rem}

Let now $S$ be a compact orientable surface, possibly disconnected and with boundary, and
$\Mod(S)$ its mapping class group.

\begin{definition} \label{def:MCG_cond}
We say that a subgroup $\Gamma < \Mod(S)$ satisfies the {\it mapping class group conditions for
primitivity} if
\begin{itemize}
\item[{\rm (i)}] it is not virtually abelian,

\item[{\rm (ii)}] it contains no finite normal subgroups, and

\item[{\rm (iii)}] there exists a surface $R$ and an embedding
$\Gamma < \Mod(R)$ such that $\Gamma$ is irreducible and acts transitively on the connected
components of $R$.
\end{itemize}
\end{definition}

\begin{theorem} \label{thm:MCG}
An infinite subgroup $\Gamma < \Mod(S)$ is primitive if and only
if it satisfies the mapping class group conditions for
primitivity.
\end{theorem}

\medskip

For convergence groups we prove the following theorem.
\begin{theorem}
\label{thm:Convergence} A countable non-elementary convergence group is primitive if and only if
it contains no finite normal subgroups.
\end{theorem}

In particular this theorem holds for subgroups of Gromov hyperbolic groups, for Kleinian groups,
or more generally for any countable group acting properly discontinuously on a Gromov hyperbolic
metric space.

\medskip

Let us recall that a group action on a tree is called {\it minimal} if there are no invariant
proper subtrees. We do not assume that the tree in the following theorem is locally finite, it
will automatically be countable though, admitting a minimal action by a countable group.

\begin{theorem}
\label{thm:Trees} Let $T$ be a tree with $|\partial T| \ge 3$. Then any countable subgroup $\Gamma
< \Aut(T)$ that acts minimally on $T$ is primitive\footnote{ Note that if $\Aut(T)$ acts minimally
on $T$ and $\partial T$ has more than two points, then $\partial T$ is actually infinite.}.
\end{theorem}

The study of maximal subgroups is related to the study of the {\it Frattini subgroup}
$\Phi(\Gamma)$ which is, by definition, the intersection of all maximal subgroups. Equivalently,
the Frattini subgroup consists of all the non generators -- elements that are expendable from any
generating set. For infinite groups, some variants of the Frattini group were introduced. These
include the near Frattini subgroup and the lower near Frattini subgroup. All of these subgroups,
including the Frattini, are contained in the subgroup $\Psi(\Gamma)$, the intersection of all
maximal subgroups of infinite index. Thus, many results are best stated in terms of the group
$\Psi(\Gamma)$.

Over the years, a few articles were published, showing that the Frattini subgroup is small for
finitely generated groups in many of the geometric settings dealt with in this article. Examples
include Platonov \cite{Platonov:Frattini} and Wehrfritz \cite{Wehrfritz:Frattini} in the linear
case, Ivanov \cite{Ivanov:Frattini},\cite[chapter 10]{Ivanov:MCG} for subgroups of mapping class
groups and Kapovich \cite{Kapovich:Frattini} for hyperbolic groups. Our analysis of maximal
subgroups gives a uniform approach to tackle all these settings simultaneously, and in the greater
generality of countable (not necessarily finitely generated) groups.

\begin{theorem}\label{thm:Frat}
Let $\Gamma$ be a countable group. If $\Gamma$ is linear in
characteristic zero, or if it is finitely generated and linear in
positive characteristic then $\Psi (\gC )$ is solvable. If $\gC$
is a subgroup of the mapping class group of a surface then
$\Psi(\Gamma)$ is solvable by finite, and $\Phi(\Gamma)$ is
solvable. If $\Gamma$ is a non-elementary convergence group then
$\Psi(\Gamma)$ is finite. Whenever $\Gamma$ acts minimally on a
tree $T$ with $|\partial T|
> 2$, the action of $\Psi(\Gamma)$ is trivial.
\end{theorem}

The part of the statement concerning minimal group actions on trees gives an answer to a question
of G. Higman and B.H. Neumann \cite{HN:Ito} on the Frattini group of an amalgamated product, in
the case of countable groups. We refer the reader to
\cite{Allenby:upper},\cite{Allenby:loc},\cite{AT:Frat_Amalgam}, \cite{Moh:Near_Frat_Amalgam} for
previous results in this vein.

\begin{corollary} (Higman-Neumann question)
\label{cor:Higman} Let $\Gamma = A*_{H}B$ be an amalgamated free product. Assume that $\gC$ is
countable and that $([A:H] -1)([B:H]-1) \ge 2$. Then $\Psi(\Gamma) < \Core_{\Gamma}(H)$. In
particular this result holds also for the Frattini, the near Frattini and the lower near Frattini
subgroups of $\Gamma$.
\end{corollary}

\subsection{Methods}

Our methods are influenced by the work of Margulis-So\u{\i}fer, who were in turn inspired by the
famous paper of Tits \cite{Tits:alternative} on the so called ``Tits alternative''. In every
finitely generated linear group $\Gamma$ which is not virtually solvable, Tits constructed a
non-abelian free subgroup. Margulis and So\u{\i}fer constructed a subgroup $\Delta$ which is free
and profinitely dense. In the current work, under a stronger assumption on the group, we construct
a subgroup $\Delta < \Gamma$ which is free and prodense, thus proving the existence of a faithful
quasiprimitive action. For finitely generated groups the existence of a faithful quasiprimitive
action is actually equivalent to the existence of a faithful primitive action (c.f.
\ref{prop:Prim=QP}). In all the cases considered in this paper, we actually prove this equivalence
also for countable non finitely generated groups.

The focus has shifted since the work of Tits. The requirement that $\Delta$ be a free group is no
longer the goal but merely a means for making sure that $\Delta \ne \Gamma$.  The specific case
where $\Gamma$ is a free group on a countable number of generators can be easily taken care of
separately (c.f. \cite{Cameron:Countable_free}). Still it is worthwhile noting that in all the
situations that we handle we actually construct a free prodense subgroup.

The proof splits into two parts, ``representation theoretical'' and ``dynamical''. The former one
consists of finding a representation in which every non-trivial normal subgroup exhibits a rich
dynamic in its action on the associated geometric boundary (e.g. the corresponding projective
space for linear group, the boundary of the tree, the limit set $\ldots$). For instance, in the
linear case we need a projective representation over some local field in which every normal
subgroup acts strongly irreducibly and the group has ``plenty'' of very proximal elements. We
borrow the relevant statement from \cite{BG:Topological_Tits}.

The latter, dynamical, part consists of constructing a free prodense subgroup $\Delta < \Gamma$.
We make sure that $\Delta$ is free by requiring the generators to satisfy the conditions of the
ping-ping Lemma \ref{lem:ping-pong1}. We add the generators one by one, making sure that at least
one generator falls in every coset of every non-trivial normal subgroup of $\Gamma$, forcing $\gD$
to be prodense. The case where $\gC$ is not finitely generated is more subtle, since we have to
construct $\gD$ inside some proper maximal subgroup of $\gC$.

In Section \ref{sec:axiom} we axiomatize the dynamical part of the proof. We state properties of a
group action on a topological space that imply primitivity. The resulting theorem is used in the
following sections to conclude the cases of convergence groups, groups acting minimally on trees,
and mapping class groups. While the argument for linear groups is conceptually similar, it is
technically more complicated in many ways. We therefore give an independent proof to the linear
case in Section \ref{sec:Linear}.

Let us note that the situation dealt here is more complicated than
the one in \cite{MS:Maximal}, as general normal subgroups are more
difficult to handle than finite index ones -- there are more of
them and they may lie deeper in the group. For instance, one
crucial point of the proof is to construct a, so called, very
proximal element in every normal subgroup. If $g\in\gC$ is some
very proximal element and $N\lhd\gC$ is a normal subgroup of
finite index then $g^{[\gC :N]}$ is a very proximal element in
$N$. If $[\gC :N]=\infty$ however, this trick does not apply
anymore. In the linear case we make a substantial use, in order to
settle this difficulty and others, of the recent work of
Breuillard and Gelander \cite{BG:Dense_Free},
\cite{BG:Topological_Tits} while for subgroups of mapping class
groups we exploit the work of Ivanov \cite{Ivanov:MCG} who proved
a Margulis-So\u{\i}fer type theorem for this setting.

\subsection{How to read this paper}
%The block diagram below describes the dependence of the chapters.

%\begin{figure}[htb]
%\begin{center}
%\resizebox{0.7\textwidth}{!}{
% \includegraphics{how_to_read4}
%}
%\end{center}
%\end{figure}%
Section \ref{sec:Linear} about linear groups can be read independently, but it might be
recommended to read first one of the easier geometric settings. A good choice would be to start
out with convergence groups reading Sections \ref{sec:first_steps}, \ref{sec:axiom},
\ref{sec:convergence}, before approaching the linear case. A reader who is interested only in the
proof of the Higman-Neumann conjecture should read Sections \ref{sec:first_steps},
\ref{sec:axiom}, \ref{sec:trees} and the relevant parts of Section \ref{sec:Frattini}, namely the
discussion in the beginning, Lemma \ref{lem:primitive_Frat} and the proof of Corollary
\ref{cor:Higman}.

\subsection{Thanks}
 We would like to thank George Glauberman for explaining to us the examples of primitive solvable
groups; Tim Riley for a discussion which lead to Corollary
\ref{cor:Stalling}; Emmanuel Breuillard and Yves de Cornulier who
pointed out errors in an earlier version of this manuscript;
Anders Karlsson who suggested that the natural setting for our
proof was that of convergence groups rather than hyperbolic
groups; R. Allenby for encouraging us to tackle the Higman-Neumann
question in a greater generality. We thank Mikl\'os Ab\'ert, Pete
Storm, Alex Furman and Gregory So\u{\i}fer for many interesting
discussions on prodense subgroups. Finally we would like to thank
the referee for carefully reading the paper and pointing out
numerus improvements and corrections.

%-------------------------------------------------------------------------------------------

\section{A strategy for constructing pro-dense
subgroups}\label{sec:strategy}\label{sec:first_steps} Constructing primitive actions amounts to
finding maximal subgroups. By Zorn's lemma any proper subgroup of a finitely generated group is
contained in a maximal one. The issue is to find a maximal subgroup which is ``not too big''.
Margulis and So\u{\i}fer construct maximal subgroups of infinite index. We are interested in
maximal subgroups which have a trivial core.

In a rather paradoxical fashion Margulis and So\u{\i}fer make sure that their subgroup is small in
one manner by requiring it to be large in a different manner. They ensure that a subgroup is of
infinite index by requiring it to be profinitely dense. The advantage of this approach is that the
property of being profinitely dense is stable when passing to bigger subgroups, and at the same
time a profinitely dense proper subgroup is always of infinite index.

Analogously, we require a group theoretic property that on one hand implies that a subgroup has a
trivial core, and on the other hand is stable under passing to bigger subgroups. This is exactly
where the prodense subgroups from Definition \ref{def:prodense} come into the picture. In fact if
$\Delta < \Gamma$ is proper and prodense then

\begin{enumerate}
\item Every subgroup containing $\Delta$ is also prodense.
\item The action $\Gamma \circlearrowleft \Gamma/\Delta$ is faithful, or equivalently,
$\Core_{\Gamma}(\Delta) = \trivgp$.
\end{enumerate}

Since, by definition, $\Delta < \Gamma$ is prodense if and only if $\Gamma \circlearrowleft \Gamma
/ \Delta$ is quasiprimitive and faithful, we have the following proposition.

\begin{proposition}\label{prop:Prim=QP}
For a finitely generated group $\Gamma$ the following three conditions are equivalent:
\begin{itemize}
\item $\Gamma$ is primitive.
\item $\Gamma$ is quasiprimitive.
\item $\Gamma$ contains a proper prodense subgroup.
\end{itemize}
\end{proposition}

The assumption that $\gC$ is finitely generated in Proposition \ref{prop:Prim=QP} can be replaced
by the weaker assumption that $\gC$ contains a finitely generated subgroup with a non-trivial
core, i.e. a finitely generated subgroup which is open in the normal topology, or by the
assumption that any prodense subgroup is contained in a maximal subgroup. However, since we
consider general countable groups we have to be more careful, and to construct a prodense subgroup
which possess the additional property that it is contained in some maximal subgroup.

In groups that have many quotients, such as hyperbolic groups, the very existence of a proper
subgroup that maps onto every proper quotient is somewhat surprising. We refer the reader to
\cite{Dixon:free_HT,MD_free_HT,AG:Generic} for constructions of prodense subgroups in free groups.
The interested reader may find stronger results for free groups in these references. In the first
two papers, highly transitive faithful actions of free groups are constructed. The third paper
constructs an action $F \circlearrowleft F/\Delta$ which is not two transitive, but for which $N
\Delta = F$ for every non-trivial subnormal subgroup $N \lhd \lhd F$. The methods in all of these
references however are very specific to free groups.

Our aim is to construct a proper pro-dense subgroup $\gD\leq\gC$ which is contained in some
maximal subgroup, whenever such a subgroup exists. Let us first explain how we construct $\gD$ to
be pro-dense, and later how we make sure that it is contained in a maximal subgroup. We want $\gD$
to project onto every proper quotient of $\Gamma$. In principle $\Gamma$ might have continuously
many normal subgroups. However, since any non-trivial normal subgroup contains the conjugacy class
of some non-trivial element, it is enough to require that $\Delta N = \Gamma$ for every normal
subgroup $N$ which is generated by a unique non-trivial conjugacy class. Since $\Gamma$ is
countable, it has only countably many conjugacy classes. In other words, unlike the profinite
topology, the normal topology is usually not countable, however, it is always second countable,
i.e. admits a countable base. Let $\mathcal{F}$ be the countable set of normal subgroups which are
generated by a non-trivial conjugacy class. Since each $N \in\mathcal{F}$ has at most countably
many cosets we can enumerate all the cosets $\{ C_{N,i}\}_{N \in \mathcal{F}, i=1 \ldots [\Gamma
:N]}$ (where $[\gC :N]\leq\aleph_0)$. We shall construct $\gD$ to be a free group on countably
many generators $\gd_{N,i}\in C_{N,i}$, one inside each coset of each $N\in\mathcal{F}$. This will
guarantee that
\begin{itemize}
\item $\gD$ maps onto every quotient $\Gamma /N,~N\in\mathcal{F}$,
and hence onto every proper quotient of $\Gamma$, and

\item $\gD \ne \Gamma$ (This is true unless $\Gamma$ happens to be
a countably generated free group, a case that can be treated separately (c.f.
\cite{Cameron:Countable_free}).),
\end{itemize}
or in other words, that $\gD$ is a proper pro-dense subgroup of $\Gamma$.

We shall make use of the following infinite variant of the classical ping-pong lemma.

\begin{lemma}\label{lem:ping-pong1}
Suppose that a group $G$ acts on a set $X$ and suppose that $g_1,g_2,\ldots$ is a sequence of
elements in $G$ such that for each $g_i$ there are four subsets $A(g_i) =A_i^+, R(g_i) = R_i^+,
A(g_i^{-1}) = A_i^-, R(g_i^{-1}) = R_i^-$ of $X$ such that the following are satisfied:
\begin{itemize}
\item $A_i^+\cap (A_i^-\cup R_i^+ )=A_i^-\cap (A_i^+\cup R_i^-)=\emptyset$,
\item $A_i^\pm\cap (A_j^\pm\cup R_j^\pm)=\emptyset,~\forall i\neq j$, and
\item $g_i\cdot (X\setminus R_i^+)\subset A_i^+,~\textnormal{and}~g_i^{-1}\cdot (X\setminus R_i^-)\subset A_i^-.$
\end{itemize}
Then the elements $\{g_i\}$ form a free generating set of a free subgroup of $G$.
\end{lemma}

We will refer to a set of elements satisfying the conditions of the lemma as a {\it ping-pong
tuple}.

In many cases, we can take $R_i^-=A_i^+$ and $R_i^+=A_i^-$. We shall then simply denote the
associated sets by $A_i,R_i$ (or $A(g_i),R(g_i)$). The three conditions in the ping-pong Lemma
\ref{lem:ping-pong1} are then replaced by the simple requirement that all the sets $A(g_i),R(g_i)$
are pairwise disjoint. For the sake of simplicity, let us restrict ourselves, for the time being,
to that situation.

We construct the free pro-dense subgroup $\gD$ in two steps.

\medskip

{\bf Step 1 (A free group intersecting every non-trivial normal subgroup).} We shall construct
$a_N\in N$ for each $N\in\mathcal{F}$ which satisfy the conditions of Lemma \ref{lem:ping-pong1},
with corresponding attracting and repelling sets ${A}(a_N),{R}(a_N)$.

\medskip

{\bf Step 2 (A free group intersecting every coset).} For each given $N\in\mathcal{F}$ we shall
construct $\gd_{N,i}\in C_{N,i}$ which exhaust the cosets $C_{N,i}$ of $N$ in $\gC$ and which
satisfy the condition of Lemma \ref{lem:ping-pong1}, and the additional requirement for the
positions of the attracting and repelling neighborhoods:
$$
 {A}(\gd_{N,i})\cup {R}(\gd_{N,i})\subset {A}(a_N).
$$
This will guaranty that the elements $\{\gd_{N,i}\}_{N\in\mathcal{F},1\leq i\leq [\Gamma
:N]\leq\aleph_0}$ satisfy the condition of Lemma \ref{lem:ping-pong1} all together.

\begin{rem}\label{srk}
Note that in the more general setup, the elements $a_N$ comes with four sets
${A}(a_N),{A}(a_N^{-1}),$ ${R}(a_N),{R}(a_N^{-1})$, and so do the $\gd_{N,i}$'s. Then one should
impose some conditions on the positions of the associated sets for $\gd_{N,i}$. The condition we
found most convenient to require in the linear case in Section \ref{sec:Linear} is:
\begin{eqnarray*}
 {A}(\gd_{N,i})\subset {A}(a_N),\ \ &~{A}(\gd_{N,i}^{-1})\subset {A}(a_N^{-1}), \\
 {R}(\gd_{N,i})\subset {R}(a_N),\ \ &~{R}(\gd_{N,i}^{-1})\subset {R}(a_N^{-1}).
\end{eqnarray*}
\end{rem}

\medskip

Finally, let us indicate how we guarantee that the prodense subgroup $\gD$ is contained in some
maximal proper subgroup: we shall construct in advance two elements $h_1,h_2$ which will be in
``general position'' with all the elements $\gd_{N,i}$. Then, after constructing the $\gd_{N,i}$,
we will add one element $\ti{c}_j$ in each non-trivial double coset $\langle h_1,h_2\rangle c_j
\langle h_1,h_2\rangle$ of the group $\langle h_1,h_2\rangle$ such that the bigger set
$\{\gd_{N,i},\ti{c}_j\}$ will still form a ping-pong tuple. Among the subgroups which contains
this big ping-pong tuple there is, by Zorn lemma, a maximal one which does not contain $\langle
h_1,h_2\rangle$. Such a subgroup must be maximal in $\gC$. Since it is maximal and prodense, the
action of $\gC$ on its cosets space is faithful and primitive.

\section{Axiomatization of the ping-pong argument}
\label{sec:axiom}

In this section, we shall formulate and prove an abstract theorem, which could be applied in
different cases as a tool to prove primitivity. In the subsequent three sections we shall apply it
to convergence groups, subgroups of mapping class groups and groups of tree automorphisms. We find
it more convenient to formulate, prove and apply, a version which is not the most general. This
modest version cannot be applied for the case of linear groups. However the proof of Theorem
\ref{thm:Linear} for linear groups is more complicated due to several other issues and we will
take care of it separately in Section \ref{sec:Linear}.

\begin{definition}
\label{def:f-prox} Let $M$ be a topological space. We call a homeomorphism $g \in \Homeo(M)$ {\it
d-contracting} if there exist disjoint open sets $A,R$ with $A \cup R \subsetneqq M$ such that $g
(M \setminus R) \subset A$. We refer to the sets $A,R$ above as {\it attracting} and {\it
repelling} open sets for $g$ respectively. Alternatively, we say that $g$ is $(A,R)$
d-contracting. Note that if $g$ is $(A,R)$ d-contracting then $g^{-1}$ is automatically $(R,A)$
d-contracting. We say that a homeomorphism $g$ is {\it f-proximal} if there are finite sets of
fixed points $g^{+} = \{a_1,\ldots,a_n\}, \ g^{-} = \{r_1,\ldots,r_n\}$ such that $g^{+} \cap
g^{-} = \emptyset$ and, for any choice of open sets $A,R$ such that $A \cap R = \emptyset$, $A
\cup R \ne M$, $A\supset g^{+}$ and $R\supset g^{-}$, there exists $m \in \N$ such that $g^{m}$ is
$(A,R)$ d-contracting. Again note that if $g$ is f-proximal then so is $g^{-1}$, with
$(g^{-1})^{+} = g^{-}, \ \ (g^{-1})^{-} = g^{+}$. Finally, we will say that a homeomorphism $g$ is
{\it proximal} if it is f-proximal and $g^{-}$ and $g^{+}$ are single points.
\end{definition}

\begin{rem} $(i)$ The notion of $d$-contracting element is similar but not
identical to the notion of contracting element defined in Section
\ref{sec:Linear}. The letter ``d'' stands for the disjointness of
the attracting and repelling neighborhoods. The notion of a
contracting element on the other hand does not require such
disjointness, but then one needs some metric on $M$ in order to
capture the contraction property.

$(ii)$ The letter $f$ in $f$-proximal stands for the finiteness of
the sets of attracting and repelling points. Usually $g^{-}$ and
$g^{+}$ will be singletons, but for mapping class groups of
disconnected surfaces we will require finitely many. Again this
notion is similar, but not identical to the notion of proximal
elements, defined in Section \ref{sec:Linear}, where a proximal
element admits a repelling hyperplane of codimension one.

$(iii)$ As long as $M$ satisfies the separation property $T_1$,
which we will always assume, it automatically follows from the
definition that $\Fix(g) = g^{-} \cup g^{+}$ for every
$f-$proximal element $g$.
\end{rem}

\begin{theorem}
\label{thm:combinatorial} Let $M$ be a regular topological space\footnote{in the sense that it is
Hausdorff, and every point can be separated from any closed set by a pair of disjoint open sets.},
$\Gamma < \Homeo(M)$ a countable group, and assume that the following hold.
\begin{enumerate}
\item \label{itm:infinite_orbits} Any orbit of any non-trivial
normal subgroup is infinite. \item \label{itm:prox} There is an f-proximal element in $\gC$. \item
\label{itm:cont->prox} For every two disjoint finite sets $S,T \subset M$ there exist open sets
$U_S,U_T$ such that $U_S \cap U_T = \emptyset$, $U_S \cup U_T \ne M$, $S \subset U_S$, $T \subset
U_T$ and every $(U_S,U_T)$ d-contracting element in $\gC$ is f-proximal.
\end{enumerate}
Then the group $\Gamma$ is primitive.
\end{theorem}

Before proceeding with the proof we shall require two lemmas.

\begin{lemma}\label{lem:finite_sets}
Assumption (\ref{itm:infinite_orbits}) of the theorem is
equivalent to the requirement: ``for every non-trivial normal
subgroup $N \lhd \Gamma$ and every two finite subsets $S, T
\subset M$ there exists an element $n \in N$ such that $nS \cap T
= \emptyset$''.
\end{lemma}

\begin{proof} Note that $\{n \in  N | nS \cap T \ne \emptyset\}$
is a finite union of cosets of stabilizers of points $N_x, \ \ x \in S$. However, a group is never
equal to a finite union of cosets of infinite index subgroups by a theorem of Neumann (see
\cite{Neumann:finitely_many_cosets}, \cite[Theorem 3.3C]{DM:Permutation_Groups}).
\end{proof}

\begin{lemma}\label{b1b2b3}
Let $M$ and $\gC$ be as in Theorem \ref{thm:combinatorial}, and
suppose that $g\in\gC$ is $(A,R)$ $d$-contracting for some $A,R
\subset M$ and $f$-proximal with $g^-\subset R,g^+\subset A$, and
$b_1,b_2,b_3\in\gC$ are three elements such that:
$$
 b_1R\cap R=b_2A\cap A=b_3^{-1}R\cap R=b_1R\cap b_3^{-1}R=\emptyset,
$$
then, for a sufficiently large $k\in\N$, the element
$$
 gb_1g^{-(1+k)}b_2g^{(k+1)}b_3g^{-1}
$$
is $(gb_1g^{-k}R,gb_3^{-1}g^{-k}R)$ $d$-contracting and $f$-proximal. Moreover,
$gb_1g^{-k}R,gb_3^{-1}g^{-k}R\subset A$.
\end{lemma}

\begin{proof}
Choose a sufficiently small neighborhood $g^{-}\subset \Omega'\subset R$ such that every
$(gb_1\Omega',gb_3^{-1}\Omega')$ $d$-contracting element is $f$-proximal. Since for $k$ large
enough $g^{-k}R\subset \Omega'$, we can actually fix such a $k$ and replace $\Omega'$ by its
subset $\Omega=g^{-k}R$. Then $R=g^k \Omega$. The element $a=gb_1g^{-(1+k)}b_2g^{(k+1)}b_3g^{-1}$
is $(gb_1\Omega,gb_3^{-1}\Omega)$ $d$-contracting since
\begin{equation} \nonumber
\begin{array}{lll}
 a(M\setminus gb_3^{-1}\Omega)&=&(gb_1g^{-k}g^{-1}b_2g)g^k(M\setminus\Omega)
     =(gb_1g^{-k}g^{-1}b_2)g(M\setminus R) \\ &\subset&
 (gb_1g^{-k})g^{-1} (b_2A)\subset gb_1g^{-k} R=gb_1\Omega.
\end{array}
\end{equation}
By the choice of $\Omega$, $a$ is also $f$-proximal. The inclusions
$gb_1\Omega,gb_3^{-1}\Omega\subset A$ are obvious.
\end{proof}

\begin{proof} [Proof of Theorem \ref{thm:combinatorial}]
We shall follow Steps 1 and 2 outlined in Section \ref{sec:strategy} and construct a prodense
subgroup $\gD$ in $\gC$, that is contained in a maximal subgroup. The reader should refer to
Figure 1 for an illustration of the big ping pong table constructed in this proof.
\begin{figure}[htb] \label{fig:pptable}
\begin{center}
\resizebox{0.75\textwidth}{!}{
 \includegraphics{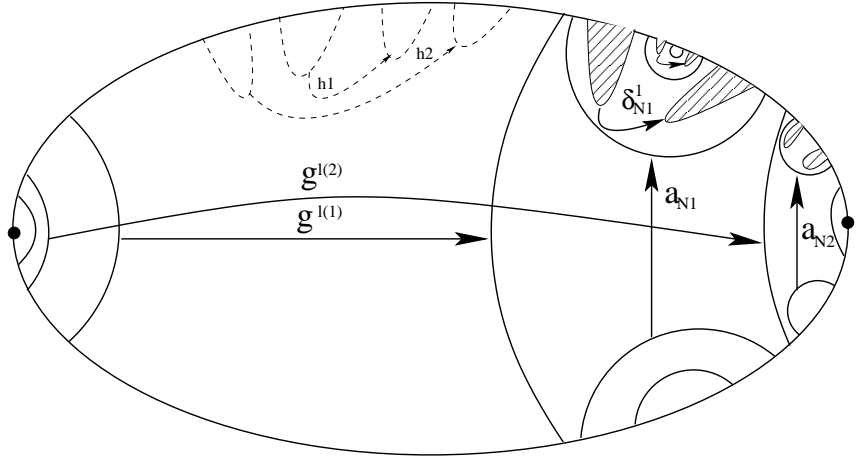}
} \caption{A big ping pong table. The arrows point from the repelling to the attracting
neighborhoods of the corresponding f-proximal elements.}
\end{center}
\end{figure}

Let $\mathcal{F} = \{N_1,N_2,\ldots\}$ be an enumeration of the set of all normal subgroups of
$\Gamma$ that are generated by a single conjugacy class. We shall artificially add two elements to
this list $N_{-1} = N_{0} = \Gamma$. We define these two artificial elements in order to guarantee
that the prodense subgroup to be constructed in Steps 1 and 2 will be contained in a maximal
proper subgroup. As $N_{-1}$ and $N_0$ play a different role than the others, we do not consider
them part of the set $\mathcal{F}$.

\medskip

\noindent {\bf Step 1} [A free set intersecting every non-trivial normal subgroup.]

By Assumption (\ref{itm:prox}), $\Gamma$ contains an $f$-proximal element. We fix such an element
$g$ once and for all, and use it throughout the argument of Step 1.

For any $m \ge -1$ we shall construct an $f$-proximal element $a_m \in N_m$ and an integer $l =
l(m)$ such that $\{a_{-1},a_0,a_1,a_2,a_3,\ldots,a_{m},g^{l}\}$ will form a ping-pong tuple. In
fact, we will require an induction hypothesis that is slightly stronger: ``there exist open
attracting and repelling neighborhoods $A_i\supset a_i^{+}$, $R_i\supset a_i^{-},~i\leq m$,
$A(g^{l(m)})\supset g^+$ and $R(g^{l(m)})\supset g^-$ with disjoint closures.'' This hypothesis is
easier to handle. At the $m$'th step\footnote{Note that for the base step of the induction $m=-2$
we can take $l(-2)=1$ and $(A(g),R(g))$ to be any pair of attracting and repelling sets for $g$.
It will follow from the proof that all the attracting and repelling neighborhoods of the elements
to be constructed are contained in $A(g)$. This fact is actually important for the argument of
Step 2.}, after $a_{-1},\ldots,a_{m-1}$ and $l(m-1),A(g^{l(m-1)}),R(g^{l(m-1)})$ are already
constructed, we shall first construct an $f$-proximal element $a_m \in N_m$ with attracting and
repelling neighborhoods satisfying $\overline{A_m},\overline{R_m} \subset A(g^{l(m-1)})\setminus
g^+$ so that it will automatically satisfy the induction hypothesis with $a_{-1},\ldots, a_{m-1}$.
Then we will pick $l(m)$ sufficiently large so that $g^{l(m)}$ will have sufficiently small
attracting and repelling set $A(g^{l(m)}),R(g^{l(m)})$, with closure disjoint from
$\overline{A_m},\overline{R_m}$.

By Lemma \ref{lem:finite_sets}, we can find $b_1,b_2,b_3\in N_m$ such that the finite sets
$(g^+\cup g^-),b_1(g^+\cup g^-),b_2(g^+\cup g^-),b_3^{-1}(g^+\cup g^-)$ are pairwise disjoint. By
regularity of the space and continuity of the action, we can find open neighborhoods $U^+\supset
g^+$ and $U^-\supset g^-$ with disjoint closures, and with $\overline{U}^+\subset A(g^{l(m-1)})$
such that the sets
$$
 b_1(\overline{U^+\cup U^-}), b_2(\overline{U^+\cup U^-}), b_3^{-1}(\overline{U^+\cup U^-}), {\overline{U^+\cup
 U^-}},
$$
are still pairwise disjoint. Since $g$ is $f$-proximal, if $j$ is sufficiently large then $g^j$ is
$(U^+,U^-)$ $d$-contracting, and by Lemma \ref{b1b2b3} we can choose $k$ such that
$$
 a_m=g^jb_1g^{-j(1+k)}b_2g^{j(k+1)}b_3g^{-j}\in N_m
$$
is $(A_m,R_m)$ $d$-contracting $f$-proximal with
$A_m=g^jb_1g^{-jk}U^-,~R_m=g^jb_3^{-1}g^{-jk}U^-$, and $A_m\cup R_m\subset U^+\subset
A(g^{l(m-1)})$. One can easily verify that $\overline{A_m}\cup \overline{R_m}$ is disjoint from
the fixed points set $g^+\cup g^-$. By regularity of $M$, we can choose a small neighborhoods
$A(g^{l(m)})\supset g^+,~R(g^{l(m)})\supset g^-$ whose closures are disjoint from
$\overline{A_m}\cup \overline{R_m}$, and then actually choose the integer $l(m)$ sufficiently
large so that $g^{l(m)}$ is $(A(g^{l(m)}),R(g^{l(m)}))$ $d$-contracting, justifying the notations
we gave to these sets. This completes the proof of the induction step, and hence the proof of Step
1.

Let us now rename the first two elements
$$
 h_1\defeq a_{-1},~h_2 \defeq a_{0},
$$
and forget about them for the time being. After finishing with Step 2, we shall use them to show
that the set constructed in that step is contained in a maximal subgroup.

\medskip

\noindent {\bf Step 2} [A free set intersecting every coset.]

Step 2 is proved by the exact same argument which proved Step 1. We shall not repeat the argument,
but only indicate the small modification that should be made in it.

As indicated in the previous section, for each given $N_m\in \mathcal{F}, \ m \ge 1$ we shall
construct cosets representatives $\gd_{N_m,i}\in C_{N_m,i}$ to be $f$-proximal elements with
attracting and repelling neighborhoods contained in $A_m$. This condition guarantee that
$\gd_{N_{m_1},i}$ and $\gd_{N_{m_2},j}$ are automatically a ping-pong pair whenever $m_1\neq m_2$.
Thus, we should only explain how to construct $\gd_{N,i}$ for each given $N\in\mathcal{F}$.

We fix $N\in\mathcal{F}$ and carry the same argument as in Step 1, almost word by word, this time
with the $f$-proximal element $a_N$, constructed in Step 1, playing the rule of $g$. Except that
here, when we choose the elements $b_1,b_2,b_3$ (in the $m$'th step of the inductive argument), we
must pick exactly one of them, say $b_1$, from the corresponding coset $C_{N,m}$, and the other
two, $b_2,b_3$, from the normal subgroup $N$ (in order that the element
$\gd_{N,m}=a_N^jb_1a_N^{-j(1+k)}b_2a_N^{j(k+1)}b_3a_N^{-j}$ will belong to the coset $C_{N,m}$).
This concludes Step 2.

\medskip

We now define
$$
 \Delta = \langle \delta_{N_m,i} \rangle_{1 \le m \le \infty, 1 \le
 i \le [\Gamma:N_m] }.
$$

This subgroup is proper and prodense in $\gC$, hence $\gC$ is quasiprimitive. In order to obtain
the primitivity of $\gC$ we shall show that $\gD$ is contained in some maximal subgroup. We shall
enlarge it to a subgroup $\ti\gD$ which can easily be seen to be contained in a maximal one.
Recall that we still have two auxiliary elements $h_1,h_2$. We shall add to the free set
$\{\gd_{N,i}\}$ countably many elements $\ti c_j$, one in each double coset $\langle
h_1,h_2\rangle c_j\langle h_1,h_2\rangle$ of $\langle h_1,h_2\rangle$, in such a way that the
larger set $\{\gd_{N,i},\ti c_j\}$ will still be free, and hence the subgroup $\ti\gD :=\langle
\gd_{N,i},\ti c_j\rangle$ will still be proper. Since $\ti\gD$ intersects any non-trivial double
coset of $\langle h_1,h_2\rangle$, an intermediate group $\ti\gD\leq A\leq\gC$ is equal to $\gC$
if and only if it contains $\langle h_1,h_2\rangle$. Since $\langle h_1,h_2\rangle$ is generated
by two elements, we can use Zorn lemma and obtain a maximal subgroup among those containing
$\ti\gD$ and not containing $\{ h_1,h_2\}$. Such a subgroup must therefore be maximal in $\gC$.

Let $\ti C_j,~j\in J$ ($\text{card}(J)\leq\aleph_0$) be the collection of all non trivial double
cosets of $\langle h_1,h_2\rangle$ in $\gC$. We shall construct $\ti c_j\in\ti C_j$ which is
$d$-contracting $f$-proximal element with
$$
 A(\ti c_j),R(\ti c_j)\subset A(h_1)\setminus h_1^+
$$
in a similar way to how the elements $a_i$ and the elements $\gd_{N,i}$ were constructed. I.e. the
attracting and repelling neighborhoods of the $\ti c_j$'s will be closer and closer but at the
same time disjoint from $h_1^+$ and altogether mutually disjoint. This will guarantee, as $A(h_1)$
is disjoint from the attracting and repelling sets of the $\gd_{N,i}$, that the large set
$\{\gd_{N,i},\ti c_j\}$ will still form a ping-pong tuple.

Denote the attracting and repelling sets of $h_i$ by $h_i^+$ and $h_i^-$ respectively. Let $c_j$
be an arbitrary element in $\ti C_j$. Multiplying $c_j$ by a sufficiently large power of $h_1$
from the left, if necessary, we may assume that $(c_j\cdot h_2^+) \cap h_2^- = \emptyset$. It
follows then that for a sufficiently large $n$, the element $h_2^n c_j h_2^n$ is $d$-contracting
with attracting and repelling neighborhoods $A(h_2^nc_jh_2^n),R(h_2^nc_jh_2^n)$ contained in those
of $h_2$. Now for every open neighborhood $U$ of $h_1^+$ we can chose $m$ sufficiently large so
that $h_1^m\cdot (A(h_2^nc_jh_2^n)\cup R(h_2^nc_jh_2^n))$ is contained in $U$. It follows that the
element $h_1^mh_2^nc_jh_2^nh_1^{-m}$ is $d$-contracting with attracting and repelling sets
contained in $U\setminus h_1^+$. Choosing $U$ small enough so that it is disjoint from the
attracting and repelling sets of $\ti c_k$ for all $k<j$, we construct, inductively, in this way
the desired double coset representative $\ti c_j$.

\end{proof}

%=================================================================================

\section{Convergence groups} \label{sec:convergence}
The theory of convergence groups, which emphasizes the dynamical-topological properties of group
actions on compacta seems as if it was tailored to accommodate our proof. We find this setting
appealing also for its wide generality. Kleinian groups, Hyperbolic groups and their subgroups,
groups acting properly discontinuously on  complete locally compact Gromov hyperbolic spaces, and
relatively hyperbolic groups, all can be realized as convergence groups. Convergence groups were
introduced by Gehring and Martin in \cite{GM:Convergence}. Gehring and Martin were studying
Kleinian groups through their action on the boundary, consequently they considered only
convergence groups acting on spheres. Later it was noticed that the definition of convergence
groups can be generalized to groups acting on general compact Hausdorff spaces. In this generality
Gromov hyperbolic groups act as convergence groups on their ideal boundary
\cite{Tukia:convergence, Tukia:Uniform_convergence_groups, Bowditch:characterisation_hyperbolic}.
More recently it was shown by Asil Yaman \cite{Yaman:Relatively_Hyp_Convergence} that relatively
hyperbolic groups can also be realized as convergence groups.

Readers who are interested in specific geometric examples, can
read the proof below considering the action of their favorite
group on the relevant boundary. The proof should make sense, with
some obvious adjustments to the terminology. Discrete subgroups of
rank one simple Lie groups are examples of convergence groups but
they are also linear groups. Hence we obtain two different, but
equivalent criteria for primitivity for these groups. We wish to
thank Anders Karlsson who suggested that the natural setting for
our proof was that of convergence groups rather than hyperbolic
groups.

The goal of this section is to prove Theorem \ref{thm:Convergence} from the introduction.

%--------------------------------------------------------------------------------------------------------

\subsection{Facts and lemmas on convergence groups}
Our survey of convergence groups follows \cite{Bowditch:Convergence}.
\begin{definition}
\label{defn:collapsing} An infinite set $\Phi$ of homeomorphisms
of a compact set $M$ is called {\it collapsing} with respect to a
pair of (not necessarily distinct) points $(a,r)$ if for every
pair of compact sets $K \subset M \setminus \{r\}$ and $L \subset
M \setminus \{a\}$, the set $\{\phi \in \Phi | \phi K \cap L \ne
\emptyset \}$ is finite. We shall then call $a$ the attracting
point and $r$ the repelling point of $\Phi$. A set $\Phi$ is {\it
collapsing} if it is a collapsing with respect to some pair of
points.
\end{definition}

\begin{definition}
\label{defn:convergence} An action of an infinite group $\Gamma$ on a topological space $M$ is
said to have the convergence property if every infinite subset $\Phi < \Gamma$ contains an
infinite subset $\Phi' \subset \Phi$ which is collapsing. A {\it convergence group} is a group
that admits a convergence action on some infinite compact Hausdorff space.
\end{definition}

It follows easily from the definition that the kernel of every convergence action is always
finite. Thus, even though we do not assume explicitly that the given convergence action is
faithful, it will automatically be almost faithful. In particular a convergence group with no
finite normal subgroups admits a faithful convergence action.

It follows from the definition of convergence group that every element of infinite order fixes
either one or two points of $M$. As a Corollary we obtain the ``usual'' classification of elements
into three mutually exclusive categories: elements of finite order are called {\it elliptic},
elements of infinite order fixing exactly one point are called {\it parabolic} and elements of
infinite order fixing two points of $M$ are called {\it loxodromic}. In dynamical terms the
loxodromic elements are proximal in the sense of Definition \ref{def:f-prox}. A similar
classification holds for subgroups.

\begin{lemma}
\label{lem:group_class} (see for example \cite{Tukia:convergence}) Let $\gC$ be a convergence
group with associated compact space $M$. Suppose an infinite subgroup $\gC' < \Gamma$ fixes some
point of $M$. Then either $\gC'$ consists entirely of elliptic and parabolic elements, or consists
entirely of elliptic and loxodromic elements. In the latter case $\Gamma$ also fixes some other
point $q \in M \setminus \{p\}$, and is virtually cyclic. Every infinite torsion group (i.e.
infinite group consisting only of elliptic elements) must fix a unique point in $M$.
\end{lemma}

\begin{definition}
\label{defn:elementary} A subgroup $\gS < \Gamma$ of a convergence group is called elementary if
it is finite or if it stabilizes a nonempty subset of $M$ with at most $2$ elements.
\end{definition}

As usual, in order to prove primitivity we have to show that every normal subgroup is big. In the
setting of convergence groups we define the {\it limit set} $L(\Gamma) \subset M$ to be the set of
all accumulation points of $\Gamma$ orbits. The limit set is also the minimal non-empty closed
$\Gamma$ invariant set. If $N \lhd \Gamma$ is an infinite normal subgroup then $L(N)$ is non empty
closed and $\Gamma$ invariant. This proves the following lemma.

\begin{lemma}
\label{lem:limitset} Let $N \lhd \Gamma$ be an infinite normal subgroup of a non-elementary
convergence group then $L(N) = L(\Gamma)$.
\end{lemma}

\begin{lemma}
\label{lem:pol->lox} If an element $g \in \Gamma$ is d-contracting, in the sense of Definition
\ref{def:f-prox}, with repelling and attracting open neighborhoods $A,R$. Then $g$ is loxodromic
and therefore proximal.
\end{lemma}

\begin{proof}\label{lem:char lox}
Both $\cap_{n \in N} g^n(M \setminus R)$ and $\cap_{n \in \N} g^{-n} (M \setminus A)$ are
non-empty closed $\langle g \rangle$ invariant sets. Thus they must be single points.
\end{proof}

The following lemma will be used only for the group $\Gamma$ itself but it is just as easy to
prove it for a general normal subgroup of $\Gamma$.

\begin{lemma}
\label{lem:lox} Let $N \lhd \Gamma$ be an infinite normal subgroup of a non-elementary convergence
group. Then $N$ contains a loxodromic element.
\end{lemma}

\begin{proof}
By Lemma \ref{lem:group_class} every infinite torsion group must fix a unique point, which is
impossible for a normal subgroup in view of Lemma \ref{lem:limitset}, hence $N$ must contain an
element of infinite order, say $n$. Let $(n^{+},n^{-})$ be the attracting and repelling fixed
points given by the Definitions (\ref{defn:convergence},\ref{defn:collapsing}) for the set $\{
n^i:i\in\N\}$. If $n^{+} \ne n^{-}$ then $n$ is loxodromic. Otherwise, by minimality of the action
of $N$ on the limit set (Lemma \ref{lem:limitset}), we can find an element $m \in N$ such that $m
n^{+} \ne n^{+}$. For large enough $j$ the element $m n^j$ will take the complement of a small
neighborhood of $n^{-}$ to a small neighborhood of $m n^{+}$, so that $m n^j$ must be loxodromic
by Lemma \ref{lem:lox}.
\end{proof}

\subsection{The proof of Theorem \ref{thm:Convergence}}
By Lemma \ref{lem:finite_normal} of the appendix, the condition is necessary.

Let $\Gamma$ be a non-elementary convergence group with no finite normal subgroups, acting
(convergently) on a compact Hausdorff space $M$. Let $L \subset M$ be the limit set. Since $L$ is
a compact Hausdorff space it is regular. We shall prove that $\gC$ is primitive by showing that
its action on $L$ satisfies the three assumptions of Theorem \ref{thm:combinatorial}. First note
that the action is faithful since the kernel is a finite normal subgroup and hence trivial. By
Lemma \ref{lem:limitset}, every normal subgroup $N \lhd \Gamma$ acts minimally on $L$ and in
particular has only infinite orbits -- this establishes Assumption (\ref{itm:infinite_orbits}).
Since, by Lemma \ref{lem:lox}, every non-elementary convergence group contains a loxodromic
element, we also have Assumption (\ref{itm:prox}). Finally, by Lemma \ref{lem:pol->lox}, every
d-contracting element of $\Gamma$ is proximal -- this gives Assumption (\ref{itm:cont->prox}).
\qed

%================================================================================

\section{Mapping class groups} \label{sec:MCG}
This section is dedicated to the proof of Theorem \ref{thm:MCG}. We take all the terminology
pertaining to mapping class groups from Ivanov's book \cite{Ivanov:MCG}. We state here a more
explicit version of Theorem \ref{thm:MCG}.

\begin{theorem} \label{thm:MCG_explicit}
Let $S$ be an orientable compact surface, $\Gamma < \Mod(S)$ an
infinite group, $\sigma = \sigma(\Gamma)$ a canonical reduction
system for $\Gamma$, $C$ a one dimensional sub-manifold in the
isotopy class of $\sigma$, $S_C=S\setminus C$ the surface $S$ cut
along $C$ and $\rho: \Mod(S) \arrow \Mod(S_C)$ the canonical
homomorphism. Let $S_C = T_1 \cup T_2 \cup \ldots \cup T_m$ be the
unique decomposition of $S_C$ as a disjoint union of subsurfaces
in such a way that $\rho(\Gamma)$ acts transitively on the set of
connected components of each $T_i$. Then the group $\Gamma$ is
primitive if and only if the following conditions are satisfied:
\begin{enumerate}
\item \label{itm:MOD_inj} The canonical map $p_i \circ \rho: \Gamma \arrow \Mod(T_i)$
is injective for some $1 \le i \le m$,
\item \label{itm:MOD_va} $\Gamma$ is not virtually abelian, and
\item \label{itm:MOD_fn} $\Gamma$ has no finite normal subgroup.
\end{enumerate}
\end{theorem}

The situation here is similar to the linear group case. The characterization of primitive
subgroups of mapping class groups of orientable compact surfaces is stated in terms of some
realization of such group as a subgroup of a mapping class group of a surface, which may not be
the original one but can be constructed from it in a few simple steps.

In our treatment of subgroups of mapping class groups we follow Ivanov \cite{Ivanov:MCG} who
proved the key theorem of the existence of pseudo-Anosov elements in infinite irreducible
subgroups of mapping class groups, and used it to prove a Margulis-So\u{\i}fer type theorem for
finitely generated subgroups of mapping class groups.

\subsection{Generalities on mapping class groups}
Let $R$ be an orientable compact surface, possibly disconnected and with boundary. We write $R =
R_1 \cup R_2 \cup \ldots \cup R_m$ as a union of connected components. Following Ivanov we
consider the action of $\Mod(R)$ on $\PMF(R)^{\#}= \PMF(R_1) \coprod \PMF(R_2) \coprod \ldots
\coprod \PMF(R_m)$, the disjoint union of the Thurston boundaries of the connected components.
Thurston's classification of homeomorphisms says that every element of the mapping class group
falls into exactly one of the following categories:
\begin{itemize}
\item \underline{periodic}, elements of finite order,
\item \underline{reducible}, preserves a one dimensional sub-manifold which is not boundary parallel, up to isotopy.
\item \underline{pseudo-anosov}, these elements exhibit f-proximal dynamics on $\PMF(R)^{\#}$.
\end{itemize}
It follows from this classification that the only mapping classes that exhibit d-contracting
dynamics on $\PMF(R)^{\#}$ are the pseudo-Anosov elements, thus we have the following corollary.

\begin{corollary} \label{cor:cont_PA}
Assume that $g \in \Mod(R)$ acts as a $d$-contracting homeomorphism on $\PMF(R)^{\#}$ then $g$ is
a pseudo-Anosov element and in particular it is f-proximal.
\end{corollary}

A similar classification holds for subgroups of the mapping class group. In particular, we call a
subgroup irreducible if it does not preserve a one dimensional sub-manifold. If a subgroup $\Gamma
< \Mod(R)$ is reducible, then one can find a realization of $\Gamma$ as a group of homeomorphisms
that actually preserve a one dimensional sub-manifold $C$. This gives rise to a canonical
homomorphism $\rho: \Gamma \arrow \Mod(R_C)$ where $R_C$ is the surface $R$ cut along $C$, and one
can check that this is well defined at the level of mapping classes. This is the reason for the
name ``reducible''-- we reduce $\Gamma$ by mapping it to a mapping class group of the ``simpler''
surface $R_C$. Note that the map $\rho$ is not always injective, but the kernel is generated by
the Dehn twists on the connected components of $C$ so it is a finitely generated abelian subgroup.

The following theorem of Ivanov is a fundamental result in the theory of mapping class groups.

\begin{theorem}(Ivanov, see \cite[Theorem 6.3 and Theorem 1 of the introduction]{Ivanov:MCG})
\label{thm:irred_PA} An infinite irreducible subgroup of the mapping class group of a surface
always contains a pseudo-Anosov element.
\end{theorem}

We shall also need a lemma.

\begin{lemma} \label{lem:Normal_irreducible}
Assume that $\Gamma < \Mod(R)$ is irreducible and acts transitively on the connected components of
$R$. Then every infinite normal subgroup $N \lhd \Gamma$ is irreducible.
\end{lemma}

\begin{proof}
This is proved in Ivanov's book in the case where $R$ is connected, see \cite[Corollary
7.13]{Ivanov:MCG}. For the general case, let $\Gamma' \lhd \Gamma$ be a finite index normal
subgroup that does not permute the connected components of $R$, and let $N' = N \cap \Gamma'$.
Since irreducibility is not sensitive to changes of finite index, $\Gamma'$ is irreducible and so
are its projections on the mapping class groups of the connected components $p_i(\Gamma') <
\Mod(R_i)$. Since $N'$ is an infinite normal subgroup, $p_i(N \cap \Gamma')$ must be infinite for
at least one coordinate $i$. However, $\Gamma$ acts transitively on the connected components and
$N' \lhd \Gamma$ so we deduce that $p_i(N)$ is infinite for every $i$. As the lemma is known for
connected surfaces, we deduce that $p_i(N')$ is irreducible for every $i$, and therefore that $N'$
and $N$ are also irreducible subgroups.
\end{proof}

\begin{corollary} \label{cor:Normal_contains_PA}
Assume that $\Gamma < \Mod(R)$ is irreducible and acts transitively on the connected components of
$R$. Then every infinite normal subgroup $N \lhd \Gamma$ contains a pseudo-Anosov element.
\end{corollary}

\begin{proof}
This is a direct consequence of Theorem \ref{thm:irred_PA} and Lemma \ref{lem:Normal_irreducible}.
\end{proof}

Let $L = L(\Gamma) = \overline{\cup\{\Fix(f):f \in \Gamma, {\textnormal{ pseudo-Anosov}}\}}
\subset \PMF^{\#}(R),$ be the canonical limit set.

\begin{corollary} \label{cor:MCG_Normal_limit}
Assume that $\Gamma < \Mod(R)$ is irreducible and acts transitively on the connected components of
$R$. Then every infinite normal subgroup $N \lhd \Gamma$ acts minimally on $L(\Gamma)$.
\end{corollary}

\begin{proof}
It is clear that $L(\Gamma)$ is contained in any closed non-empty $\Gamma$-invariant subset of
$\PMF^{\#}(S)$. From Lemma \ref{lem:Normal_irreducible} it follows that $L(N) \ne \emptyset$ for
every infinite normal subgroup $N \lhd \Gamma$. Since $N$ is normal, $L(N)$ is a closed
$\Gamma$-invariant subset of $\PMF^{\#}(R)$. Thus $L(N) = L(\Gamma)$. Clearly $N$ acts minimally
on $L(N)$.
\end{proof}

\subsection{Necessary conditions} Let $\Gamma < \Mod(S)$ be a quasiprimitive group.
Conditions \ref{itm:MOD_va} and \ref{itm:MOD_fn} of Theorem \ref{thm:MCG_explicit} follow
immediately from Proposition \ref{prop:MN} and Lemma \ref{lem:finite_normal}. Since $\Gamma$ is a
subgroup of a mapping class group it is residually finite (see for example \cite[Exercise
11.1]{Ivanov:MCG}) and by Corollary \ref{cor:RF_non_banal} it cannot be banal. In other words if
$N,M \lhd \Gamma$ are non-trivial normal subgroups then $[M,N] \ne \trivgp$. Let $\rho: \Gamma
\arrow \Mod(S_C)$ be the canonical homomorphism. The kernel $\ker (\rho )$ is generated by the
Dehn twists along the components of $C$, and in particular it is abelian. Since $\Gamma$ is not
banal $\ker \rho = \trivgp$, and $\rho: \Gamma \arrow \Mod(S_C)$ is injective. Now write $S_C =
T_1 \cup T_2 \cup \ldots \cup T_m$ as a disjoint union of surfaces where $\Gamma$ acts
transitively on the connected components of $T_i$ for each $i$. This gives rise to an embedding
$\Gamma < \Mod(T_1) \times \Mod(T_2) \times \ldots \times \Mod(T_m)$. Since $\Gamma$ is not banal
it maps injectively into one of these factors. This verifies Condition \ref{itm:MOD_inj} and
completes the proof of the necessary conditions of Theorem \ref{thm:MCG_explicit}.

\subsection{Sufficient conditions}
Assume that a countable group $\Gamma < \Mod(R)$ is irreducible, not virtually cyclic, contains no
finite normal subgroups and acts transitively on the connected components $R_1,\ldots,R_m$ of $R$.
We will prove that $\Gamma$ is primitive by showing that all the conditions of Theorem
\ref{thm:combinatorial} hold for the action $\Gamma \circlearrowleft L(\Gamma)$ where $L =
L(\Gamma) = \overline{\cup_{f \in \Gamma, {\textnormal{ pseudo-Anosov }}}\Fix(f)} \subset
\PMF^{\#},$ is the canonical limit set. Note that since $\Gamma$ is not virtually abelian and
irreducible, the limit set $L$ is infinite. Furthermore, the action of $\Gamma$ on $L$ is faithful
since $\Gamma$ contains no finite normal subgroups.

Since $L$ is a compact Hausdorff space it is regular. By Corollary \ref{cor:cont_PA} Assumption
(\ref{itm:cont->prox}) of Theorem \ref{thm:combinatorial} holds. By Corollary
\ref{cor:MCG_Normal_limit}, every normal subgroup $N \lhd \Gamma$ acts minimally on $L$ and
therefore Assumption (\ref{itm:infinite_orbits}) of Theorem \ref{thm:combinatorial} also holds.
Finally, Assumption (\ref{itm:prox}) of Theorem \ref{thm:combinatorial} holds by Corollary
\ref{cor:Normal_contains_PA}. Thus all the conditions of Theorem \ref{thm:combinatorial} hold and
the group $\Gamma$ is primitive. This concludes the proof of Theorem \ref{thm:MCG_explicit}.

%================================================================================

\section{Groups acting on trees} \label{sec:trees}
\subsection{Generalities about trees}
A very similar analysis can be carried out for group actions on trees. Let $T$ be a tree which is
locally finite or locally countable, $\Aut(T)$ its automorphism group and $\partial T$ the
boundary:
$$
 \partial T = \{f: \N \arrow T| {\textrm{ f is an infinite geodesic ray}}\}/\sim
 $$
where two infinite rays are equivalent if their images eventually coincide
$$
 f \sim g \Leftrightarrow \exists m,n {\textnormal{ such that }} f(m+i) = g(n+i),~ \forall i \ge 0.
$$
A unique geodesic path $[x,y]$ connects any two points $x,y \in T \cup \partial T$. This path
might be finite or infinite on any of the sides depending on whether the points are in the tree or
on the boundary. If $x \ne y \in T$ are two vertices, we define the shadow
$$
 \Shadow_{x \arrow y} = \{ \eta \in \partial T | y \in [x,\eta]\}.
$$
The collection of all shadows forms a basis of open neighborhoods for a topology on $\partial T$.
When $T$ is locally finite, $T \cup \partial T$ is its natural compactification, but in general
$\partial T$ is not compact. In any case, $T$ embeds as a dense open discrete subset into $T \cup
\partial T$ with the natural topology, and the action of $\Aut(T)$ on $T$ extends canonically to a
continuous action on $T \cup \partial T$.

If $\Gamma < \Aut(T)$ is any subgroup, we define the {\it limit set} $L(\Gamma) \subset
\partial T$ as the set of all accumulation points of orbits of $\Gamma$ on $T$. The
limit set is the minimal $\Gamma$-invariant closed subset of $\partial T$.

\begin{definition}
\label{defn:minimal_group} A subgroup $\Gamma < \Aut(T)$ is called {\it minimal} if it admits no
invariant subtree.
\end{definition}

If $\Gamma$ is minimal then $L(\Gamma) = \partial T$ and $\Gamma$ acts on $\partial T$ minimally
since the convex core of the limit set is always an invariant subtree. Conversely, assuming that
$\Aut(T)$ acts minimally on $T$, a subgroup $\Gamma < \Aut(T)$ is minimal if and only if it acts
minimally on the boundary and does not fix a vertex or a geometric edge (this follows easily from
\cite[Propositions 7.1 and 7.5]{Bass:Covering}).

\begin{lemma}
\label{lem:tree_minimal} Let $\Gamma < \Aut(T)$ be a minimal group and $\trivgp \ne N \lhd
\Gamma$. Then $N$ contains a hyperbolic element, $L(N) = \partial T$ and $N$ acts minimally on
$\partial T$.
\end{lemma}

\begin{proof}
Assume that $N$ does not contain any hyperbolic element then by
standard results for group actions on trees (see
\cite{Serre:Trees}) $N$ has to fix a point $x \in T$ (this point
might be a vertex or the center of a geometric edge). Since $N$ is
normal it must fix pointwise the closed convex hull of the orbit
$\overline{\Conv(\Gamma \cdot x)}$, which is everything by
minimality of the $\Gamma$ action. Since the action of $\Gamma$ is
faithful, it follows that $N = \trivgp$.

When $N \ne \trivgp$, it contains a hyperbolic element and therefore $L(N) \ne \emptyset$. Since
$N \lhd \Gamma$, the limit set $L(N)$ is closed and $\Gamma$-invariant so $L(N) = L(\Gamma) =
\partial T$. Finally $N$ acts minimally on $\partial T$ because any group acts minimally on its limit set.
\end{proof}

In the study of actions on trees that are not locally finite, we encounter for the first time the
situation where a d-contracting element does not have to be proximal. For locally finite trees
this does not occur because every elliptic element preserves some natural measure on the boundary.
However in the locally countable case, elliptic elements can be d-contracting. In fact, one can
visualize a contraction behavior even on the set of nearest neighbors of the fixed vertex. The
following lemma shows that the more delicate Assumption (\ref{itm:cont->prox}) of Theorem
\ref{thm:combinatorial} still holds.

\begin{lemma}
\label{lem:tree_pol->lox} Let $A,R \subset \partial T$ be open sets. Suppose that these sets are
``far away from each other'' in the sense that there exist a path in the tree $x = x_0,x_1,x_2
\ldots x_n = z$ and $\eta \in \partial T$ satisfying the following conditions:
\begin{itemize}
\item $A \subset \Shadow_{x \arrow z}$ and $R \subset \Shadow_{z \arrow x}$.

\item $[\eta,x] \cap [\eta,z]\ni x_j$ for some $0 < j <n$.
\end{itemize}
Let $g \in \Aut(T)$ be $(A,R)$ contracting in the sense that $g(M \setminus R) \subset A$. Then
$\partial T\setminus (A\cup R)\neq\emptyset$, and the element $g$ is hyperbolic, and hence
proximal.
\end{lemma}

\begin{proof}
Clearly, $\eta \in \partial T\setminus (A\cup R)$. In order to prove that $g$ is hyperbolic, it is
enough, by \cite[Lemma 6.8]{Bass:Covering}, to exhibit one edge $(s,t)$ that is not fixed by $g$
and which is coherent with $(g\cdot s,g\cdot t)$, in the sense that the path $[s,g\cdot s]$
contains exactly one of the vertices $t$ or $g\cdot t$. In our case one can take the edge
$(x_0,x_1)$.
\end{proof}

\subsection{The proof of Theorem \ref{thm:Trees}.}
All we have to do is verify that the conditions of Theorem \ref{thm:combinatorial} hold for the
action $\Gamma \circlearrowleft \partial T$. The space $\partial T$ is a regular because it admits
a basis of clopen sets (actually it is easy to see that it is even metrizable). Assumptions
(\ref{itm:infinite_orbits}) and (\ref{itm:prox}) of Theorem \ref{thm:combinatorial} follow from
Lemma \ref{lem:tree_minimal}. Assumption (\ref{itm:cont->prox}) follows from Lemma
\ref{lem:tree_pol->lox}. \qed

%================================================================================

\section{Linear Groups}
\label{sec:Linear} This section is devoted to the proof of the following theorem.
\begin{theorem}  \label{thm:Linear_nonbanal}
A countable non-torsion linear group $\gC$ which is not banal is primitive if and only if it
satisfies the linear conditions for primitivity \ref{def:lin_cond}.
\end{theorem}
The proof of Theorem \ref{thm:Linear} is then concluded in the appendix where it is shown that
banal groups are either of affine or of diagonal type and that these groups admit a unique
quasiprimitive action, which is actually primitive.

\subsection{Necessary conditions}
\label{subsec:Necessary} Let $\Gamma$ be a countable infinite quasiprimitive linear group. If
$\Gamma$ is finitely generated then it is not banal by Corollary \ref{cor:RF_non_banal}. Otherwise
$\Gamma$ is not banal by assumption. By Lemma \ref{lem:finite_normal} $\Gamma$ contains no finite
normal subgroups. To complete the proof of the necessary conditions we will show that if $\Gamma$
is a linear group with no finite normal subgroups, and in which no two non-trivial normal
subgroups commute elementwise, then it satisfies the linear conditions for primitivity given in
Definition \ref{def:lin_cond}.

Let $f: \Gamma \arrow \GL_n(k)$ be a faithful linear representation over some algebraically closed
field. We will show how, in a few simple steps, we can modify the representation $f$ to get a new
representation which satisfies the conditions of Theorem \ref{thm:Linear}.

Let $\G = \overline{f(\Gamma )}^Z$ be the Zariski closure and $\G^{\circ}$ the connected component
of the identity in $\G$. The intersection of $f(\Gamma)$ with the solvable radical of $\G^{\circ}$
has to be trivial. Indeed if $\Gamma$ had a non-trivial solvable subgroup $S \lhd \Gamma$, then it
would also have a non-trivial abelian normal subgroup -- the last non-trivial group in the derived
series of $S$, and in particular it would be banal. Dividing by the solvable radical we obtain a
new faithful representation $f_1$ such that the Zariski closure of the image $\G_1 =
\overline{f_1(\Gamma)}^{Z}$ is semisimple. The lack of normal abelian subgroups also implies that
$f_1(\Gamma) \cap Z(\G_1) = \trivgp$. Take $f_2=\Ad_{\G_1}\circ f_1$, then $f_2$ is still
faithful, and $\G_2 = \overline{f_2(\Gamma)}^{Z}$ is center-free and semisimple. Write
$\G_2^\circ=\prod_{i=1}^k \H_i$ according to the partition determined by the orbits of the factors
of $\G_2^\circ$ under the $\G_2$ action by conjugations, and denote by $\H_i$ the images of $\G_2$
in $\Aut(\H_i^\circ)$. Then $\G_2$ and hence $\gC$ embeds into the direct product $\prod_{i=1}^k
\H_i$, where for each $\H_i$, the connected component $\H_i^\circ$ is a direct product of
isomorphic simple algebraic groups and the action of $\H_i$ (as well as of $\G_2$ and of $\gC$) on
$\H_i^\circ$ permutes the simple factors transitively. Since $\Gamma$ does not contain elementwise
commuting normal subgroups we can divide by the product of all the $\H_i,i \neq i_0$ for some $1
\leq i_0 \leq k$ and obtain a faithful representation $f_3$ of $\Gamma$ into $\H = \H_{i_0}$. Now
since $\H^\circ$ is center-free and has finite index in $\H$, the kernel of the action of $\Gamma$
on $\H^\circ$ by conjugation composed with $f_3$ is a finite normal subgroup of $\Gamma$ and hence
trivial. We thus obtain a faithful representation of $\Gamma$ into the linear algebraic group
$\Aut(\H^\circ )$ which satisfies the linear conditions for primitivity from Definition
\ref{def:lin_cond}.

\medskip

We shall now aim at proving the sufficient conditions part of Theorem \ref{thm:Linear}. In order
to simplify, we shall restrict ourselves throughout the argument to the case where $\gC$ is
finitely generated. At the end of this section, we shall indicate the changes needed to be made
when $\gC$ is not assumed to be finitely generated.

%================================================================================

\subsection{Some preliminaries about projective transformations over local fields}

In this paragraph we shall review some definitions and results from \cite{BG:Dense_Free} and
\cite{BG:Topological_Tits} regarding the dynamical properties of projective transformations which
we shall use in the proof.

Let $k$ be a local field and $\left\| \cdot\right\| $ the standard norm on $k^{n}$, i.e. the
standard Euclidean norm if $k$ is Archimedean and $\left\| x\right\| =\max_{1\leq i\leq n}|x_{i}|$
where $x=\sum x_{i}e_{i}$ when $k$ is non-Archimedean and $(e_{1},\ldots ,e_{n})$ is the canonical
basis of $k^{n}$. This norm extends in the usual
way to $\Lambda ^{2}k^{n}$. Then we define the \textit{standard metric} on $%
\mathbb{P}(k^{n})$ by $d([v],[w])=\frac{\left\| v\wedge w\right\|
}{\left\| v\right\| \left\| w\right\|}$, where $[v]$ denotes the
projective point corresponding to $v\in k^n$. Unless otherwise
specified all our notation will refer to this metric, for example
$B_{\nu}(v)$ will denote the ball of radius $\nu$ around a point
$v \in \mathbb{P}(k^{n})$. With respect to this metric, every
projective transformation is bi-Lipschitz on $\mathbb{P}(k^{n})$.
For $\epsilon \in (0,1)$, we call a projective transformation
$g\in $PGL$_{n}(k)$ {\it $\epsilon $-contracting} if there exist a
point $v_{g}\in {\mathbb{P}}^{n-1}(k),$ called an attracting point
of $g,$ and a projective hyperplane $H_{g}$, called a repelling
hyperplane of $g$, such that $g$ maps the complement of the $\epsilon $%
-neighborhood of $H_{g}\subset \mathbb{P}(k^{n})$ (the repelling neighborhood of $g$) into the
$\epsilon $-ball around $v_{g}$ (the attracting neighborhood of $g$). We say that $g$ is $\epsilon
${\it-very contracting} if both $g$ and $g^{-1}$ are $\epsilon $-contracting. A
projective transformation $g\in $PGL$_{n}(k)$ is called {\it $%
(r,\epsilon )$-proximal} ($r>2\epsilon >0$) if it is $\epsilon $-contracting with respect to some
attracting point $v_{g}\in \mathbb{P}(k^{n})$ and some repelling hyperplane $H_{g}$, such that
$d(v_{g},H_{g})\geq r$. The transformation $g$ is called {\it $(r,\epsilon )$-very proximal} if
both $g$ and $g^{-1}$ are $(r,\epsilon )$-proximal. Finally, $g$ is
simply called {\it proximal} (resp. {\it very proximal}) if it is $(r$%
{\it $,\epsilon )$}-proximal (resp. $(r${\it$,\epsilon )$}-very proximal) for some
$r>2${\it$\epsilon >0$. }

The attracting point $v_{g}$ and repelling hyperplane $H_{g}$ of an $%
\epsilon $-contracting transformation are not uniquely defined. Yet, if $g$ is proximal we have
the following nice choice of $v_{g}$ and $H_{g}$.

\begin{lemma}\label{fix}(Lemma 3.2 of \cite{BG:Topological_Tits})
Let $\epsilon \in (0,\frac{1}{4})$. There exist two constants $%
c_{1},c_{2}\geq 1$ (depending only on the local field $k$) such that if $g$ is an $(r,\epsilon
)$-proximal transformation with $r\geq c_{1}\epsilon $ then it must fix a unique point
$\overline{v}_{g}$ inside its attracting neighborhood and a unique projective hyperplane
$\overline{H}_{g}$ lying inside its repelling neighborhood\footnote{by this we mean that if $v,H$
are any couple of a pointed a hyperplane with $d(v,H)\geq r$ s.t. the completion of the
$\gep$-neighborhood of $H$ is mapped under $g$ into the $\gep$-ball around $v$, then
$\overline{v}_g$ lies inside the $\gep$-ball around $v$ and $\overline{H}_g$ lies inside the
$\gep$-neighborhood around $H$}.
Moreover, if $r\geq c_{1}\epsilon ^{2/3}$%
, then the positive powers $g^{n}$, $n\geq 1$, are $(r-2\epsilon ,(c_{2}\epsilon
)^{\frac{n}{3}})$-proximal transformations with respect to these same $\overline{v}_{g}$ and
$\overline{H}_{g}$.
\end{lemma}

In what follows, whenever we add the article \textit{the} (or {\it the canonical}) to an
attracting point and repelling hyperplane of a proximal transformation $g$, we shall
mean these fixed point $\overline{v}_{g}$ and fixed hyperplane $\overline{H}%
_{g}$ obtained in Lemma \ref{fix}. Moreover, when $r$ and $\gep$ are given, we shall denote by
$A(g),R(g)$ the $\gep$-neighborhoods of $\overline{v}_g,\overline{H}_g$ respectively. In some
cases, we shall specify different attracting and repelling sets for a proximal element $g$. In
such a case we shall denote them by $\mathcal{A}(g),\mathcal{R}(g)$ respectively. This means that
$$g \big(\PP(k^n)\setminus\mathcal{R}(g)\big)\subset\mathcal{A}(g). $$
If $g$ is very proximal and we say that
$\mathcal{A}(g),\mathcal{R}(g),\mathcal{A}(g^{-1}),\mathcal{R}(g^{-1})$ are specified attracting
and repelling sets for $g,g^{-1}$ then we shall always require additionally that
$$
 \mathcal{A}(g)\cap\big(\mathcal{R}(g)\cup\mathcal{A}(g^{-1})\big)=
 \mathcal{A}(g^{-1})\cap\big(\mathcal{R}(g^{-1})\cup\mathcal{A}(g)\big)=\emptyset.
$$

\medskip

Using proximal elements, one constructs free groups with the following variant of the classical
ping-pong lemma.

\begin{lemma}\label{lem:ping-pong}
Suppose that $\{ g_i\}_{i\in I}\subset\PGL_n(k)$ is a set of very proximal elements, each
associated with some given attracting and repelling sets for itself and for its inverse. Suppose
that for any $i\neq j,~i,j\in I$ the attracting set of $g_i$ (resp. of $g_i^{-1}$) is disjoint
from both the attracting and repelling sets of both $g_j$ and $g_j^{-1}$, then the $g_i$'s form a
free set, i.e. they are free generators of a free group.
\end{lemma}

A set of elements which satisfy the condition of Lemma \ref{lem:ping-pong} with respect to some
given attracting and repelling sets will be said to form a ping-pong set (or a ping-pong tuple).

\medskip

Given a contracting element, one can construct a proximal one using the following lemma (c.f.
\cite[Section 3 and Proposition 3.8]{BG:Dense_Free}).

\begin{lemma}\label{lem:contracting->very-proximal}
Suppose that $G\leq \PGL_n(k)$ is a group which acts strongly irreducibly (i.e. does not stabilize
any finite union of projective hyperplanes) on the projective space $\PP (k^n)$. Then there are
constants
$$
 \gep (G),r(G),c(G)>0
$$
such that if $g\in G$ is an $\gep$-contracting transformation for some $\gep <\gep (G)$ then for
some $f_1,f_2\in G$ the element $gf_1g^{-1}f_2$ is $(r(G),c(G)\gep)$-very proximal.
\end{lemma}

The following characterization of contracting elements is proved in \cite[Proposition 3.3, and
Lemmas 3.4, 3.5]{BG:Dense_Free}.

\begin{lemma}\label{lem:contracting=Lipschitz}
There exists some constant $c$, depending only on $k$,
 such that for any $\gep\in (0,\frac{1}{4})$ and $d\in(0,1)$,
\begin{itemize}
\item if $g\in\PGL_n(k)$ is $(r,\gep)$-proximal for $r>c_1\gep$,
then it is $c\frac{\gep^2}{d^2}$-Lipschitz outside the $d$-neighborhood of the repelling
hyperplane, and vice versa \item if $g$ is $\gep^2$-Lipschitz on some open neighborhood then it is
$c\gep$-contracting.
\end{itemize}
Here $c_1$ is the constant given by Lemma \ref{fix}.
\end{lemma}

The main ingredient in the method we use for generating free
subgroups is a projective representation whose image contains
contracting elements and acts strongly irreducibly. The following
theorem is a particular case of Theorem 4.3 from
\cite{BG:Topological_Tits}. Note that a similar statement appeared
also earlier in \cite{MS:Maximal}.

%The following theorem is a weak version of Theorem 4.3 from \cite{BG:Topological_Tits}.

\begin{theorem}\label{thm:good-representation}
Let $K$ be a field and $\mathbb{H}$ an algebraic $K$-group for which the connected component
$\mathbb{H}^{\circ}$ is not solvable, and let $\Gamma <\mathbb{H}$ be a Zariski dense finitely
generated subgroup. Then we can find a number $r>0$, a local field $k$, an embedding
$K\hookrightarrow k$, an integer $n$, and a strongly irreducible projective representation $\rho
:\mathbb{H}(k)\rightarrow PGL_{n}(k)$ defined over $k$, such that for any $\epsilon \in
(0,\frac{r}{2})$ there is $g\in\Gamma\cap\mathbb{H}^{\circ}$ for which $\rho (g)$ acts as an
$(r,\epsilon )$-very proximal transformation on $\mathbb{P}(k^{n})$.
\end{theorem}

\subsection{Sufficient conditions}

Let $K$ be an arbitrary field and $\Gamma\leq\GL_n(K)$ a finitely generated group for which the
connected component $\G^\circ$ of $\G=\overline{\Gamma}^Z$ is a direct product of simple $K$
algebraic groups and the action of $\Gamma$ on $\G^\circ$ by conjugation is faithful and permutes
the simple factors of $\G^\circ$ transitively. One important property which follows immediately
from this condition is that for any non-trivial normal subgroup $N\lhd\Gamma$ we have
$\G^\circ\subset\overline{N}^Z$.

We shall prove that $\Gamma$ is primitive by constructing a pro-dense subgroup $\gD <\Gamma$
following the same guiding lines described in Sections \ref{sec:strategy} and \ref{sec:axiom}.
Theorem \ref{thm:good-representation} supplies us with some local field $k$, and a strongly
irreducible algebraic projective representation $\rho :\Gamma\to \PGL_n(k)$ such that $\rho (g)$
is $(r,\gep_0)$-very proximal for some $g\in \Gamma$ with $\gep_0<(\frac{r}{c_1})^{3/2}$ i.e.
$\rho (g)$ satisfies the conditions of Lemma \ref{fix}. We shall fix this $g$ and use it
throughout the proof. Denote by $\overline{v}_g,\overline{H}_g,
\overline{v}_{g^{-1}},\overline{H}_{g^{-1}}$ the attracting points and repelling hyperplane of $g$
and $g^{-1}$ respectively.

Let $\mathcal{F}$ be the countable set of normal subgroups $N\lhd\Gamma$ which are generated by a
non-trivial conjugacy class, and for each $N\in\mathcal{F}$ let $\{ C_{N,i}\}_{i=1}^{[\Gamma :N]}$
denote the cosets of $N$ in $\gC$. We shall find one element $\gd_{N,i}\in C_{N,i}$ for each
$C_{N,i}$ such that the $\gd_{N,i}$ will satisfy the condition of Lemma \ref{lem:ping-pong}, and
hence will form a free set, and take
$$
 \gD=\langle\gd_{N,i}:N\in\mathcal{F},i\in\{1,2,\ldots ,[\Gamma :N]\}\rangle.
$$
We shall do that in two steps.

\medskip

{\bf Step 1.} We shall construct $a_N\in N$ for each $N\in\mathcal{F}$ with $\rho (a_N)$
satisfying the conditions of Lemma \ref{lem:ping-pong}, i.e. they will be very proximal elements
with  attracting and repelling sets $\mathcal{A}(a_N), \mathcal{A}(a_N^{-1}), \mathcal{R}(a_N),
\mathcal{R}(a_N^{-1})$, such that each attracting set $\mathcal{A}(\cdot )$ is disjoint from the
union of the attracting and repelling sets of all $a_{N'},~N'\in\mathcal{F}\setminus\{ N\}$.

\medskip

{\bf Step 2.} For each given $N\in\mathcal{F}$ we shall construct $\gd_{N,i}\in C_{N,i}$ which
exhaust the cosets $C_{N,i}$ of $N$ and which satisfy the condition of Lemma \ref{lem:ping-pong},
and the additional requirement for the positions of the attracting and repelling neighborhoods:
\begin{eqnarray*}
 \mathcal{A}(\gd_{N,i})\subset\mathcal{A}(a_N),\ \ &~\mathcal{A}(\gd_{N,i}^{-1})\subset\mathcal{A}(a_N^{-1}),\\
 \mathcal{R}(\gd_{N,i})\subset\mathcal{R}(a_N),\ \ &~\mathcal{R}(\gd_{N,i}^{-1})\subset\mathcal{R}(a_N^{-1}).
\end{eqnarray*}

This will guarantee that the elements $\{\gd_{N,i}\}_{N\in\mathcal{F},1\leq i\leq [\Gamma
:N]\leq\aleph_0}$ satisfy the condition of Lemma \ref{lem:ping-pong} all together.

%=======================================================================================

\subsection{Step 1}

For the argument of Step 1, we shall number the elements of the countable (or finite) set
$\mathcal{F}$ by $N_1,N_2,\ldots$. We shall construct the elements $a_j=a_{N_j}$ recursively with
respect to some specified attracting and repelling neighborhoods
$\mathcal{A}(a_j),\mathcal{R}(a_j)$ and $\mathcal{A}(a_j^{-1}),\mathcal{R}(a_j^{-1})$ for $a_j$
and $a_j^{-1}$, and find some $\nu_j>0$, and $n_j\in \N$ such that $\rho(g^{n_j})$ is very
proximal with respect to the $\nu_j$-neighborhoods of the attracting points
$\overline{v}_g,\overline{v}_{g^{-1}}$ and the repelling hyperplanes
$\overline{H}_g,\overline{H}_{g^{-1}}$, and $\{a_1,\ldots ,a_j,g^{n_j}\}$ form a ping-pong tuple
with respect to the specified attracting and repelling neighborhoods.

\medskip

{\it In the $j^{\textnormal{th}}$ step we shall construct $a_j$
such that the contracting sets for $a_j$  (resp. $a_j^{-1}$) will
be contained in the $\nu_{j-1}$ neighborhood of $\overline{v}_g$
(resp. $\overline{v}_{g^{-1}}$) but disjoint from
$\overline{H}_{g^{-1}}$ (resp. $\overline{H}_g$). Similarly the
repelling neighborhoods of $a_j$ (resp. $a_j^{-1}$) will be
contained in the $\nu_{j-1}$-neighborhood of $\overline{H}_{g}$
(resp. $\overline{H}_{g^{-1}}$) while being disjoint from
$\overline{v}_{g^{-1}}$ (resp. $\overline{v}_{g}$). Then choose
$\nu_j$ small enough so that the $\nu_j$-neighborhoods of
$\overline{v}_g,\overline{v}_{g^{-1}},\overline{H}_g,\overline{H}_{g^{-1}}$
are disjoint from the specified neighborhoods of all $a_i, i\leq
j$ and choose large enough $n_j$ so that $g^{n_j}$ is very
proximal with those $\nu_j$-neighborhoods.}

%{\it In the $j^{\textnormal{th}}$ step we shall construct $a_j$ with contracting and repelling
%sets which are contained in those of $g^{n_{j-1}}$ (i.e. in the $\nu_{j-1}$ neighborhoods of
%$\overline{v}_g,\overline{v}_{g^{-1}},\overline{H}_g,\overline{H}_{g^{-1}}$) but disjoint from the
%centers (i.e. from $\overline{v}_g,\overline{v}_{g^{-1}},\overline{H}_g,\overline{H}_{g^{-1}}$),
%and then choose $\nu_j$ small enough so that the $\nu_j$-neighborhoods of
%$\overline{v}_g,\overline{v}_{g^{-1}},\overline{H}_g,\overline{H}_{g^{-1}}$ are disjoint from the
%specified neighborhoods of all $a_i, i\leq j$ and choose large enough $n_j$ so that $g^{n_j}$ is
%very proximal with those $\nu_j$-neighborhoods.}

\medskip

\begin{figure}[htb]
\begin{center}
\resizebox{0.5\textwidth}{!}{
 \includegraphics{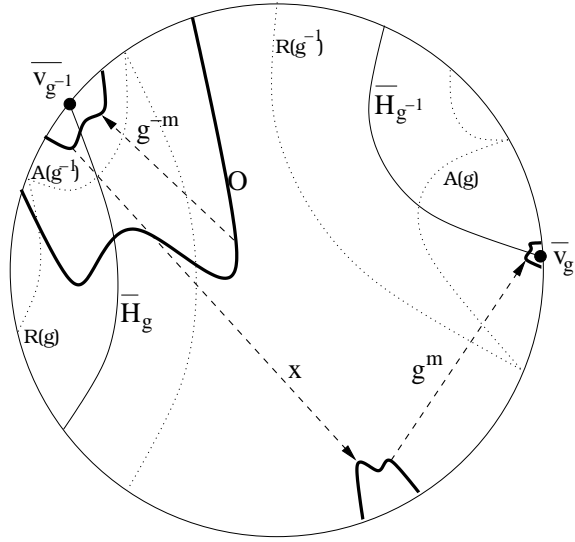}
} \caption{Constructing a contracting element, by verifying a Lipschitz condition on an open
set.}\label{fig:contracting}
\end{center}
\end{figure}

As $N_j$ is normal in $\Gamma$, it follows that
$(\overline{\Gamma}^Z)^\circ\subset\overline{N_j}^Z$ and hence $\rho (N_j)$ acts strongly
irreducibly on $\PP (k^n)$. In particular we can find $x_j\in N_j$ such that $\rho (x_j)$ moves
the attracting point $\overline{v}_{g^{-1}}$ of $g^{-1}$ outside the repelling hyperplane
$\overline{H}_g$ of $g$. Now consider the element
$$
 y_j=g^{m_j}x_jg^{-m_j}.
$$
We claim that if $m_j$ is sufficiently large then $\rho (y_j)$ is $\gep_j$-Lipschitz and hence
$c\sqrt{\gep_j}$-contracting, for arbitrarily small $\gep_j>0$. Indeed, it follows from Lemma
\ref{fix} that for a large enough $m_j$ the element $\rho (g^{-m_j})$ is $\gep_j$-very proximal
with the same attracting point and repelling hyperplane as $\rho (g^{-1})$. Then by Lemma
\ref{lem:contracting=Lipschitz} $\rho (g^{-m_j})$ is $C_1\gep_j^2$-Lipschitz on some small open
neighborhood $O$ of $\overline{v}_{g^{-1}}$, for some constant $C_1$ depending on $\rho(g)$. If we
take $O$ to be small enough and $m_j$ large enough then $\rho (x_jg^{-m_j})(O)$ has some positive
distance from the repelling $\gep_j$-neighborhood (around $\overline{H}_g$) of $\rho (g^{m_j})$.
By Lemma \ref{lem:contracting=Lipschitz} again, the element $\rho (g^{m_j})$ is
$C_2\gep_j^2$-Lipschitz on $\rho (x_jg^{-m_j})(O)$. Now since $\rho (x_j)$ is bi-Lipschitz with
some constant $C_3$ depending on $\rho (x_j)$, it follows that $\rho (y_j)=\rho
(g^{m_j}x_jg^{-m_j})$ is $C_4\gep_j^4$-Lipschitz on $O$, where $C_4$ depends only on $\rho (g)$
and $\rho (x_j )$. If we require also $\gep_j\leq\frac{1}{C_4}^{\frac{1}{3}}$ then $\rho (y_j)$ is
$\gep_j$-Lipschitz on $O$ as we claimed. Then by the second part of Lemma
\ref{lem:contracting=Lipschitz}, it follows that $\rho (y_j)$ is $c\sqrt{\gep_j}$-contracting.
This argument is illustrated in Figure \ref{fig:contracting}.

Since $\rho (N_j)$ acts strongly irreducibly, we can apply Lemma
\ref{lem:contracting->very-proximal}, assuming that $\gep_j$ is sufficiently small (i.e. that
$c\sqrt{\gep_j}\leq\gep (N_j)$) and obtain a very-proximal element $p_j\in N_j$. We shall assume
further that $\gep_j$ is small enough so that $p_j$ satisfies the condition of Lemma \ref{fix}
(for this we have to require that $\frac{r(N_j)}{[c\gep_j^{\frac{1}{2}}c(N_j)]^{\frac{2}{3}}}\geq
c_1$).

Since $N_j$ acts strongly irreducibly on $\PP (k^n)$, we can replace $p_j$ by some conjugate of it
and assume that $\rho (p_j)$ and $\rho (g)$ are in general position, i.e.
\begin{eqnarray*}
 \{\overline{v}_g,\overline{v}_{g^{-1}}\}\cap\big( \overline{H}_{p_j}\cup \overline{H}_{p_j^{-1}}\big)& = &
 \{\overline{v}_{p_j},\overline{v}_{p_j^{-1}}\}\cap\big( \overline{H}_g\cup \overline{H}_{g^{-1}}) \\
 & = &
 \{\overline{v}_g,\overline{v}_{g^{-1}}\}\cap\{\overline{v}_{p_j},\overline{v}_{p_j^{-1}}\}=
 \emptyset.
\end{eqnarray*}
Then, taking $\nu_j'\leq\nu_{j-1}$ small enough, and $l_j$ sufficiently large, we may assume that
the powers $p_j^{l_j},g^{l_j}$ form a ping-pong pair with respect to the $\nu_j'$-neighborhoods of
their canonical attracting points and repelling hyperplanes. We may also assume that $l_j$ is
large enough so that $g^{-(l_j+n_{j-1})}\big( B_{\nu_{j-1}}(\overline{v}_{g^{-1}})\big)\subset
B_{\nu_j'}(\overline{v}_{g^{-1}})$. Set
$$
 a_j=g^{l_j}p_j^{l_j}g^{-l_j}.
$$
Then $a_j\in N_j$ since $N_j$ is normal.

Now if $A(p_j^{l_j}),A(p_j^{-l_j}),R(p_j^{l_j}),R(p_j^{-l_j})$ are the attracting and repelling
neighborhoods for $p_j^{l_j}$ then
\begin{eqnarray*}
 \mathcal{A}(a_j)=g^{l_j}A(p_j^{l_j}), \ \ & \mathcal{A}(a_j^{-1})=g^{l_j}A(p_j^{-l_j}), \\
 \mathcal{R}(a_j)=g^{l_j}R(p_j^{l_j}), \ \ & \mathcal{R}(a_j^{-1})=g^{l_j}R(p_j^{-l_j})
\end{eqnarray*}
form attracting and repelling sets for the very proximal element $a_j$. Moreover,
$\mathcal{A}(a_j)\cup\mathcal{A}(a_j^{-1})=g^{l_j}\big( A(p_j^{l_j})\cup A(p_j^{-l_j})\big)
\subset B_{\nu_j'}(\overline{v}_g)$ and in particular they are disjoint from $\cup_{i=1}^{j-1}
\Big( \mathcal{A}(a_i)\cup \mathcal{A}(a_i^{-1})\cup \mathcal{R}(a_i)\cup
\mathcal{R}(a_i^{-1})\Big)$. On the other hand, since $\nu_j'\leq\nu_{j-1}$ for any $i<j$, we have
$g^{-l_j}\big(\mathcal{A}(a_i)\cup \mathcal{A}(a_i^{-1})\big) =
g^{-l_j+n_{j-1}}g^{n_{j-1}}\big(\mathcal{A}(a_i^{-1})\cup\mathcal{A}(a_i)\big)\subset B_{\nu_j'}(
\overline{v}_{g^{-1}})$ and hence disjoint from $\mathcal{R}(p_j^{l_j})\cup
\mathcal{R}(p_j^{-l_j})$. Thus $\cup_{i=1}^{j-1}\big(\mathcal{A}(a_i)\cup
\mathcal{A}(a_i^{-1})\big)$ is disjoint from $\mathcal{R}(a_j)\cup \mathcal{R}(a_j^{-1})=g^{l_j}
\left( R(p_j^{l_j}) \cup R(p_j^{-l_j}) \right)$.

Finally, since $p_j$ is in general position with $g$, also
$a_j=g^{l_j}p_j^{l_j}g^{-l_j}$ is. Hence, replacing $a_j$ with
some large power of it, if needed, we may assume that $a_j$ and
some large power $g^{n_j}$ of $g$ form a ping-pong pair. We may
also assume that $n_j\geq n_{j-1}$ and the attracting and
repelling neighborhoods of $g^{n_j}$ are the $\nu_j$-neighborhoods
of its canonical attracting points and repelling hyperplanes for
some positive number $\nu_j\leq\nu_j'$. Therefore $\{ a_1,\ldots
a_{j-1},a_j,g^{n_j}\}$ form a ping-pong tuple.

This concludes the proof of the recursive argument, and hence of Step 1. We shall drop now the
numbering $N_i$ for the elements of $\mathcal{F}$ and the notation $a_i=a(N_i)$. In the sequel, we
shall use the notation $a_N$ for the proximal element associated to $N\in\mathcal{F}$.

%==============================================================================

\subsection{Step 2}

We shall now fix $N\in\mathcal{F}$ and construct for each coset $C_{N,i}$ of $N$ in $\Gamma$ an
element $\gd_{N,i}$ and associated repelling and attracting sets $\mathcal{A}(\gd_{N,i})$,
$\mathcal{A}(\gd_{N,i}^{-1})$, $\mathcal{R}(\gd_{N,i})$, $\mathcal{R}(\gd_{N,i}^{-1})$ such that
\begin{enumerate}
\item $\rho (\gd_{N,i})$ acts as a very proximal element with
respect to these neighborhoods, i.e.
\begin{eqnarray*}
 \rho (\gd_{N,i})\big(\PP (k^n)\setminus\mathcal{R}(\gd_{N,i})\big) & \subset &\mathcal{A}(\gd_{N,i}) ~\textnormal{and}~ \\
 \rho (\gd_{N,i}^{-1})\big(\PP (k^n)\setminus\mathcal{R}(\gd_{N,i}^{-1})\big)
 & \subset & \mathcal{A}(\gd_{N,i}^{-1}).
\end{eqnarray*}
\item The elements $\gd_{N,i},~1\leq i\leq [\Gamma :N]$, form a
ping-pong tuple with respect to the corresponding specified
attracting and repelling sets. \item All the attracting and
repelling sets are contained in the corresponding sets of $a_N$,
i.e.
\begin{eqnarray*}
 \mathcal{A}(\gd_{N,i}) \subset \mathcal{A}(a_N), \ \ & \mathcal{A}(\gd_{N,i}^{-1})\subset\mathcal{A}(a_N^{-1}), \\
 \mathcal{R}(\gd_{N,i}) \subset \mathcal{R}(a_N), \ \ & \mathcal{R}(\gd_{N,i}^{-1})\subset\mathcal{R}(a_N^{-1}).
\end{eqnarray*}
\end{enumerate}
Throughout this argument we shall use the very proximal element $a_N$ in a similar (but different)
way to the use of $g$ in the proof of Step 1.

We are first going to construct a sequence $(\gb_i)_{i=1}^{[\gC :N]}$ of proximal elements in $N$
which satisfy the conditions (1),(2),(3) above. After that we shall replace them by elements
$\gd_{N,i}$ which will form a set of cosets representatives (i.e. $\gd_{N,i}\in C_{N,i}$) and
which have almost the same dynamics as the $\gb_i$'s, and in particular satisfy (1),(2),(3) above.

Since $N\cap\G^{\circ}$ is Zariski dense in $\mathbb{G}^{\circ }$, it acts strongly irreducibly on
$\PP (k^n)$ and we may pick an element $\gamma\in N$ such that
\begin{equation} \nonumber
\left\{
\begin{array}{ll}
 \rho (\gamma)\overline{v}_{{a_N}}, & ~\rho (\gamma)\overline{v}_{{a_N}^{-1}}, \\
 \rho(\gamma^{-1})\overline{v}_{a_N}, & ~\rho (\gamma^{-1})\overline{v}_{{a_N}^{-1}}
\end{array}
 \right\}
 \cap %
 \big(\overline{H}_{{a_N}}\cup\overline{H}_{{a_N}^{-1}}\cup\{\overline{v}_{{a_N}},%
\overline{v}_{{a_N}^{-1}}\}\big)=\emptyset.
\end{equation}

Now consider the element $\ga_{m_1}={a_N}^{m_1}\gamma {a_N}^{m_1}$. When ${m_1}$ is large enough,
$\ga_{m_1}$ acts on ${\mathbb{P}} (k^n)$ under $\rho$ as a very proximal transformation, whose
repelling neighborhoods lie inside the repelling neighborhoods of $a_N$ and whose attracting
neighborhoods lie inside the attracting neighborhoods of ${a_N}$. We can certainly assume that
$\rho (\ga_{m_1})$ satisfies the conditions of Lemma \ref{fix}. Hence $\ga_{m_1}$ fixes some
attracting points $\overline{v}_{\ga_{m_1}},\overline{v}_{\ga_{m_1}^{-1}}$ which are close to, but
distinct from $\overline{v}_{a_N},\overline{v}_{{a_N}^{-1}}$
respectively. Similarly the repelling neighborhoods of $\ga_{m_1},%
\ga_{m_1}^{-1}$ lie inside the repelling neighborhoods of $a_N,a_N^{-1}$, and the repelling
hyperplanes $\overline{H}_{\ga_{m_1}},\overline{H}_{\ga_{m_1}^{-1}}$ are close to that of $a_N$.
We claim that for all large enough $m_1$:
\begin{eqnarray}  \label{999}
\nonumber \{\overline{v}_{a_N},\overline{v}_{{a_N}^{-1}}\}\cap\big(\overline{H}%
_{\ga_{m_1}}\cup\overline{H}_{\ga_{m_1}^{-1}}\big)=\emptyset, \\
 \{\overline{v}_{\ga_{m_1}},\overline{v}_{\ga_{m_1}^{-1}}\}\cap%
\big(\overline{H}_{a_N}\cup\overline{H}_{{a_N}^{-1}}\big)=\emptyset.
\end{eqnarray}

Let us explain, for example, why $\overline{v}_{{a_N}^{-1}}\notin\overline{H}%
_{\ga_{m_1}}$ and why $\overline{v}_{\ga_{m_1}}\notin\overline{H}%
_{{a_N}^{-1}}$ (the other six conditions are similarly verified). Apply $%
\ga_{m_1}$ to the point $\overline{v}_{{a_N}^{-1}}$. As ${a_N}$ stabilizes $%
\overline{v}_{{a_N}^{-1}}$ we see that
\begin{equation} \nonumber
\ga_{m_1}(\overline{v}_{{a_N}^{-1}})={a_N}^{m_1}\gamma {a_N}^{m_1}(\overline{v}%
_{{a_N}^{-1}})={a_N}^{m_1}\gamma (\overline{v}_{{a_N}^{-1}}).
\end{equation}
Now, by our assumption, $\gamma (\overline{v}_{{a_N}^{-1}})\notin\overline{H}_{a_N}$.
Moreover when ${m_1}$ is large, ${a_N}^{m_1}$ is $\epsilon_{m_1}$%
-contracting with $\overline{H}_{{a_N}^{m_1}}=\overline{H}_{a_N},~\overline{v}%
_{{a_N}^{m_1}}=\overline{v}_{a_N}$ and $\epsilon_{m_1}$ arbitrarily small. Hence, we
may assume that $\gamma (\overline{v}_{{a_N}^{-1}})$ is outside the $\epsilon_{m_1}$-repelling neighborhood of ${a_N}^{m_1}$. Hence $\ga_{m_1}(%
\overline{v}_{{a_N}^{-1}})={a_N}^{m_1}\big(\gamma
(\overline{v}_{{a_N}^{-1}})\big)$ lie near $\overline{v}_{a_N}$
which is far from $\overline{H}_{\ga_{m_1}}$ hence
$\ga_{m_1}\overline{v}_{a_N^{-1}}\notin\overline{H}_{\ga_m}$.
Since $\overline{H}_{\ga_{m_1}}$ is invariant under $\ga_{m_1}$,
we conclude that
$\overline{v}_{{a_N}^{-1}}\notin\overline{H}_{\ga_{m_1}}$.

To show that $\overline{v}_{\ga_{m_1}}\notin\overline{H}_{{a_N}^{-1}}$ we shall apply
${a_N}^{-2{m_1}}$ to $\overline{v}_{\ga_{m_1}}$. If ${m_1}$ is
very large then $\overline{v}_{\ga_{m_1}}$ is very close to $\overline{v}%
_{a_N}$, and hence also ${a_N}^{m_1}(\overline{v}_{\ga_{m_1}})$ is very close to $%
\overline{v}_{a_N}$. As we assume that $\gamma$ takes $\overline{v}_{a_N}$ outside $%
\overline{H}_{{a_N}^{-1}}$, we get (by taking ${m_1}$ sufficiently large) that $%
\gamma$ also takes ${a_N}^{m_1}\overline{v}_{\ga_{m_1}}$ outside $\overline{H}%
_{{a_N}^{-1}}$. Taking ${m_1}$ even larger if necessary we get that ${a_N}^{-{m_1}}$ takes $\gamma
{a_N}^{m_1}\overline{v}_{\ga_{m_1}}$ to a small neighborhood of $\overline{v}_{{a_N}^{-1}}$. Hence
\begin{equation} \nonumber
{a_N}^{-2{m_1}}\overline{v}_{\ga_{m_1}}={a_N}^{-2{m_1}}\ga_{m_1}\overline{v}%
_{\ga_{m_1}}={a_N}^{-{m_1}}\gamma {a_N}^{m_1}\overline{v}_{\ga_{m_1}}
\end{equation}
lies near $\overline{v}_{{a_N}^{-1}}$. Since $\overline{H}_{{a_N}^{-1}}$ is ${a_N}^{-2{%
m_1}}$ invariant and is far from $\overline{v}_{{a_N}^{-1}}$, we conclude that $%
\overline{v}_{\ga_{m_1}}\notin\overline{H}_{{a_N}^{-1}}$.

Now it follows from (\ref{999}) and Lemma \ref{fix} that for every $%
\epsilon_1>0$ we can take $j_1$ sufficiently large so that ${a_N}^{j_1}$ and $%
\ga_{m_1}^{j_1}$ are $\epsilon_1$-very proximal transformations, and the $%
\epsilon_1$-repelling neighborhoods of each of them are disjoint from the $%
\epsilon_1$-attracting points of the other, and hence they form a ping-pong pair. Set
$\gb_1=\ga_{m_1}^{j_1}$.

In a second step, we construct $\gb_2$ in an analogous way to the first step, working with
${a_N}^{j_1}$ instead of ${a_N}$. In this way we would get $\gb_{2}$ which is $\epsilon_2$-very
proximal, and play ping-pong with $a_N^{j_1j_2}$. Moreover, by construction, the
$\epsilon_2$-repelling neighborhoods of $\gb_2$
lie inside the $\epsilon_1$-repelling neighborhoods of $a_N^{j_1}$, and the $%
\epsilon_2$-attracting neighborhoods of $\gb_2$ lie inside the $\epsilon_1$%
-attracting neighborhoods of $a_N^{j_1}$. Hence the three elements $\gb_1,~\gb_2$ and
$a_N^{j_1j_2}$ form a ping-pong 3-tuple.

We continue recursively and construct the desired sequence $(\gb_{j})$.

Let now $x_i\in C_{N,i}$ be arbitrary coset representatives. Since $N$ is normal, cosets are
identified with double cosets, and we can multiply the $x_i$'s from both sides by elements of $N$.

Now since $N$ acts strongly irreducibly on $\PP(k^n)$ we can multiply $x_{j}$ on the left and on
the right by some elements of $N$ so that, if we call this new element $x_{j}$ again, $\rho
(x_{j})\overline{v}_{\gb_{j}}\notin \overline{H}_{\gb_{j}}$ and $\rho
(x_{j}^{-1})\overline{v}_{\gb_{j}^{-1}}\notin \overline{H}_{\gb_{j}^{-1}}$.

Finally set
$$
 \gd_{N,j}=\gb_{j}^{l_{j}}x_{j}\gb_{j}^{l_{j}}
$$
for some positive power $l_{j}$. Then $\gd_{N,j}\in Nx_{j}N=x_jN=C_{N,j}$
Moreover, if we take $l_{j}$ large enough, it will act on $%
\mathbb{P}(V_{K})$ as a very proximal transformation whose attracting and repelling neighborhoods
are contained in those of $\gb_{j}$. Therefore, the $\gd_{N,j} $'s also form an infinite ping-pong
tuple. This concludes the proof of Step 2 and hence of Theorem \ref{thm:Linear} in the finitely
generated case. \qed

\medskip

Consider now the case where $\gC$ is not finitely generated. In fact our argument for constructing
a free prodense subgroup works almost word by word. The main change is in Theorem
\ref{thm:good-representation} from \cite{BG:Topological_Tits} which remains true under our
assumption on the Zariski closure of $\gC$ and that $\gC$ is not a torsion group, if we replace
the local field $k$ by a valuation field. The reason is that under our assumption the following
lemma holds.

\begin{lemma}
$\gC$ admits a Zariski dense finitely generated subgroup.
\end{lemma}

This Lemma is trivial in the 0 characteristic case. In the positive characteristic case it follows
by an argument similar to the proof of Lemma 5.6 in \cite{BG:Topological_Tits}.

We refer the reader to the beginning of Section 6 in \cite{BG:Topological_Tits} for a more
detailed explanation of a similar situation where non-finitely generated groups are considered.

\medskip

Finally let us also note that in order to prove primitivity rather
than just quasiprimitivity we can argue in the exact same way as
we did in Section \ref{sec:axiom}, i.e. add to $\mathcal{F}$ two
artificial elements $N_{-1}=N_0=\gC$ and start the inductive
argument in Step 1 with $m=-1$. The free prodense subgroup we
construct is then guaranteed to be contained in a maximal one.

%======================================================================================================================

%======================================================================================================================

\section{Frattini subgroups} \label{sec:Frattini}
\subsection{Generalities on Frattini and Frattini-like subgroups}
\begin{definition}
The {\it Frattini} subgroup $\Phi(\Gamma)$ of a group $\Gamma$ is the intersection of all proper
maximal subgroups of $\Gamma$. If no proper maximal subgroups exist then we define $\Phi(\Gamma) =
\Gamma$. Equivalently one can define the Frattini group as the group consisting of all
non-generators, where a {\it non-generator} is an element that is expendable from any set of
generators, i.e. $\langle A,\gamma \rangle = \Gamma~~\Rightarrow~~\langle A \rangle = \Gamma.$
Similarly we define $\Psi(\Gamma)$ to be the intersection of all maximal subgroups of infinite
index. The {\it near Frattini group} $\mu(\Gamma)$ is the intersection of near-maximal subgroups,
namely subgroups of infinite index in $\Gamma$ which are not contained in any other subgroup of
infinite index. The {\it lower near Frattini group} $\lambda(\Gamma)$ is defined as the group of
all near non-generators -- an element $\gamma \in \Gamma$ is called a {\it near non-generator} if
$[\Gamma : \langle A,\gamma \rangle] < \infty \Rightarrow [\Gamma: \langle A \rangle] < \infty \ \
\forall A \subset \Gamma$.
\end{definition}

Clearly we have inclusions of the form $\Phi(\Gamma) \leq \Psi(\Gamma)$ and $\lambda(\Gamma) \leq
\mu(\Gamma) \leq \Psi(\Gamma)$. If $f: \Gamma \arrow \Sigma$ is a surjective homomorphism then the
pull back of a maximal subgroup or of a subgroup of infinite index still retains the same
property. This immediately implies that $f(\Phi(\Gamma)) \subset \Phi(\Sigma)$, and similar
inclusions for all the other Frattini-like subgroups. All of these Frattini-like subgroups are
characteristic, and in particular normal.

One of the first observations of Frattini was that $\Phi(\Gamma)$ is nilpotent when $\Gamma$ is
finite. In fact the same argument, that came to be known as the Frattini argument, proves the
following lemma.
\begin{lemma}
\label{lem:finite_Frattini} If $\Phi(\Gamma)$ is finite, then it is nilpotent.
\end{lemma}
%\noindent The following theorem is more difficult but of similar nature.
%\begin{theorem}
%\label{thm:fg_solvable} If $\Phi(\Gamma)$ is finitely generated and solvable then it is nilpotent.
%\end{theorem}

Many works are dedicated to proving that in various geometric settings the Frattini subgroup, or
other Frattini-like subgroups are nilpotent, or otherwise small. Platonov and independently
Wehrfritz prove that $\Phi(\Gamma)$ is nilpotent when $\Gamma$ is a finitely generated linear
group \cite{Platonov:Frattini}, \cite{Wehrfritz:Frattini}. Platonov actually treats the more
general case where the matrix coefficients of $\Gamma$ are contained in some finitely generated
subring of the field. Ivanov proves the same theorem for finitely generated subgroups of mapping
class groups \cite{Ivanov:Frattini},\cite[Chapter 10]{Ivanov:MCG}. Kapovich
\cite{Kapovich:Frattini} proves that $\Phi(\Gamma)$ is finite (and therefore nilpotent) if
$\Gamma$ is a finitely generated subgroup of a word hyperbolic group. Many results were proved for
groups acting on trees and especially for amalgamated free products, see for example
\cite{Allenby:loc},\cite{Allenby:upper},\cite{Moh:Near_Frat_Amalgam},\cite{AT:Frat_Amalgam} and
the references therein.

The classification of primitive groups yields a unified approach to the proof of these and many
other results of the same nature. Moreover the primitive group approach allows us to dispose of
the finite generation assumption in all of the above settings. The key observation is the
following:

\begin{lemma}
\label{lem:primitive_Frat} $\Psi(\Gamma) = \trivgp$ for every primitive group $\Gamma$.
\end{lemma}

\begin{proof}
$\Psi(\Gamma)$ is a normal subgroup of $\Gamma$ which is contained in every maximal subgroup. For
a primitive group we can point out a specific maximal subgroup that does not contain any
non-trivial normal subgroup -- a maximal prodense subgroup.
\end{proof}

We go ahead and prove Theorem \ref{thm:Frat}. We will break the statement into a few small
theorems.

\begin{theorem}
\label{thm:conv_Frat} $\Psi(\Gamma)$ is finite, for any countable non-elementary convergence group
$\Gamma$.
\end{theorem}

\begin{proof}
Being a non-elementary convergence group, $\gC$ has a maximal finite normal subgroup $N$, and
$\gC/N$ is a non-elementary convergence group without non-trivial finite normal subgroups. By
Theorem \ref{thm:Convergence}, $\Gamma / N$ is primitive. Let $f: \Gamma \arrow \Gamma/N$ the
quotient map. By Lemma \ref{lem:primitive_Frat}, $f(\Psi(\Gamma))\leq\Psi(\Gamma/N) = \trivgp$.
Therefore $\Psi(\Gamma) < N$ and hence finite.
\end{proof}

This generalizes Kapovich's theorem in three ways. First it treats convergence groups rather than
subgroups of hyperbolic groups, second, it covers the case of countable groups rather than
finitely generated ones, and third, it establishes the finiteness of $\Psi(\Gamma)$ rather than
that of $\Phi(\Gamma)$. As $\Phi(\Gamma) \leq \Psi(\Gamma)$ the Frattini subgroup will also be
finite, and therefore automatically nilpotent by Lemma \ref{lem:finite_Frattini}.

\begin{theorem}\label{thm:lin-fra}
Let $\gC\leq\GL_n(k)$ be a countable linear group. In case
$\text{char}(K)>0$ assume further that $\gC$ is finitely
generated\footnote{Actually a weaker assumption is needed here,
namely that $\gC$ modulo the radical is non-torsion.}. Then
$\Psi(\Gamma)$ is solvable.
\end{theorem}

\begin{proof}
Let $\Gamma < \GL_n(k)$ be as above. Without loss of generality we
take $k$ to be algebraically closed. Let $G =
\overline{\Gamma}^{Z}$ be the Zariski closure. Let $R$ be the
solvable radical of $G$, $p: G \arrow H:=G/R$ the quotient map,
and $\gC_1=p(\gC )\leq H$. As $\Psi (\gC )\subset p^{-1}(\Psi
(\gC_1))$ it is enough to show that $\Psi (\gC_1)$ is trivial.

Write $H = H_1 \times H_2 \times \ldots \times H_m$ of groups
$H_i$ which satisfy the condition of Theorem \ref{thm:Linear},
i.e. the connected component $H_i^\circ$ is a direct product of
isomorphic simple center-free groups and the adjoint action of
$H_i$ on $H_i^\circ$ is faithful and permutes the simple factors
of $H_i^\circ$ transitively. Denote by $p_i: H \arrow H_i$ the
projection on the $i^{th}$ factor. By Theorem \ref{thm:Linear}
$p_i(\Gamma_1)$ is primitive and by Lemma \ref{lem:primitive_Frat}
$\Psi(p_i(\Gamma_1)) = \trivgp \quad \forall 1 \le i \le m$. Since
$p_i(\Psi(\Gamma_1)) \leq \Psi (p(\gC_1))$ we conclude that
$\Psi(\Gamma_1) = \trivgp$.
\end{proof}

Unlike Platonov's and Wehrfritz's theorems, Theorem
\ref{thm:lin-fra} treats countable groups which need not be
finitely generated, or satisfy any conditions on the matrix
coefficients (at least in characteristic 0). It also establishes
the statement for the bigger group $\Psi(\Gamma)$. Note that
Platonov and Wehrfritz prove that $\Phi(\Gamma)$ is nilpotent.
This however is no longer true once we leave the realm of finitely
generated groups. The following example is given by Philip Hall in
\cite[page 327]{Hall:Frattini}. Let $p,q$ be two primes such that
$q \equiv 1 \mod p^2$, let $C = \{z \in \C \ | \ \exists n {\text{
such that }} z^{q^n} = 1\}$ be the Pr\"ufer $q$-group, and set
$\omega$ to be a primitive $(p^2)^{th}$-root of unity in $\Z_{q}$.
Then the cyclic group $\langle a \rangle = \Z/p^2\Z$ acts on $C$
by $c \arrow c^\omega$ and the semidirect product $\Gamma =
\Z/p^2\Z \ltimes C$ is a linear group with a unique maximal
subgroup $\Phi(\Gamma) = \Z/p\Z \ltimes C$ which is not even
locally nilpotent. In fact since $\omega^p - 1$ is a unit in
$\Z_q$ then $[C_n,a^p] = C_n$ where $C_n < C$ is the unique
subgroup of order $q^n$.

Finally the statement about trees in Theorem \ref{thm:Frat} is a direct consequence of Theorem
\ref{thm:Trees} and Lemma \ref{lem:primitive_Frat}.

Let us finish this section by showing how one concludes Corollary \ref{cor:Higman}.

\begin{proof}[Proof of Corollary \ref{cor:Higman}]
Consider the amalgamated product $G = A*_{H}B$, and let $T$ be the
corresponding Bass Serre tree. The action of $G$ on $T$ is
minimal. The assumption $([A:H] -1)([B:H]-1) \ge 2$ is a short way
of saying that both of these indices are at least $2$ and that one
of them is at least $3$ and the Bass Serre tree is neither finite
nor an infinite line. Thus $|\partial T| = \infty > 3$, and we can
apply Theorem \ref{thm:Frat} and conclude that $\Psi(G)$ acts
trivially on $T$. In particular $\Psi(G) \leq H$, and since
$\Psi(G) \lhd G$ we conclude that $\Psi(G) \leq \Core_{G}(H)$.
\end{proof}

Note that the exact same proof holds for groups splitting as an HNN extension. Thus we can deduce
the following corollary.

\begin{corollary}\label{cor:Stalling}
Every finitely generated group with more than one end has a finite nilpotent Frattini subgroup.
\end{corollary}

\begin{proof}
A group with two ends is virtually cyclic. By Stalling's theorem if a group has more than two ends
(and hence infinitely many ends) it splits over a finite group. Thus the Frattini subgroup is
finite and by Lemma \ref{lem:finite_Frattini} it is also nilpotent.
\end{proof}

\begin{remark}
Many similar conclusions can be deduced. For example it is known that if a group has a finitely
generated solvable Frattini then the Frattini is in fact nilpotent (see for example
\cite{Robinson:Book}). Thus if a countable group splits over a finitely generated solvable group
then it has a nilpotent Frattini.
\end{remark}

Finally let us note that all the results proved in this Section can be improved.

\begin{rem}
The proof of Lemma \ref{lem:primitive_Frat} actually shows that
there exists some maximal subgroup with trivial core, we
therefore conclude that:
\begin{itemize}
\item a countable non-elementary convergence group has a maximal
subgroup with a finite core,

\item a countable linear group (of characteristic 0) has a
subgroup which is a finite intersection of maximal subgroups whose
core is solvable,

\item a countable non-trivial amalgamated product $A*_{H}B$ has a
maximal subgroup whose core is contained in $H$,
\end{itemize}
and so on...
\end{rem}

%=========================================================================================

\section{Further questions} \label{sec:questions}
A prodense subgroup $\Delta < \Gamma$ is by definition a group
that maps onto every proper quotient of $\Gamma$, so intuitively
it should be easier to find prodense subgroups in groups that do
not have many normal subgroups. E.g. any subgroup (even the
trivial one) of a simple group is prodense. A more instructive
example is $\PSL_2(\Z [1/p])$. It satisfies Margulis' normal
subgroup theorem as well as the congruence subgroup property, so
its family of quotients is limited to finite congruence quotients,
and indeed, its finitely generated subgroup $\PSL_2(\Z)$ is
maximal and prodense (which is the same as maximal core-free).
Hyperbolic groups on the other hand have many factor groups. The
following question is natural.

\begin{question} \label{conj:fg_pro-d}
Is every prodense subgroup of a non-elementary hyperbolic group
finitely generated.
\end{question}

The simplest examples of non-elementary hyperbolic groups are free groups and surface groups, i.e.
lattices in $\SL_2(\R)$. These groups have the LERF property -- every finitely generated subgroup
is closed in the profinite topology (see \cite{Scott:LERF} for surface groups). So in free groups
and in surface groups a finitely generated subgroup can not even be profinitely dense. There are
examples of hyperbolic groups that do not have the LERF property (see \cite{Kapovich:Howson}).

The next natural candidates are hyperbolic 3-manifold groups.
Question \ref{conj:fg_pro-d} has recently been answered
affirmatively in this case independently by Minasyan and by
Glasner-Souto-Storm.

\begin{theorem} \label{thm:SL(2,C)} (\cite{Minasyan:qc}, \cite{GSS:Lattices})
Let $\Gamma < \SL_2(\C)$ be a lattice, and $\Delta < \Gamma$ a maximal subgroup of infinite index,
or a prodense subgroup. Then $\Delta$ is not finitely generated.
\end{theorem}

Note that Theorem \ref{thm:SL(2,C)} holds also for non-uniform
lattices in $\SL_2(\C)$ which are not Gromov hyperbolic groups,
but still have many quotients. In fact any non-uniform lattice in
$\SL_2(\C)$ can be mapped onto many uniform lattices in
$\SL_2(\C)$. On the other hand, the statement about maximal
subgroups is special to lattices in $\SL_2(\C)$ and does not hold
in general for hyperbolic groups. Indeed by a famous theorem of
Rips \cite{Rips:fp_sequence}, every finitely presented group
$\Sigma$ can be placed in a short exact sequence:
$$1 \arrow N \arrow \Gamma \arrow \Sigma \arrow 1$$
where $\Gamma$ is a hyperbolic group and $N$ is generated by two
elements, as an abstract group. By choosing $\Sigma$ to be any
finitely presented group with a finitely presented maximal
subgroup of infinite index (e.g. $\PSL_2(\Z[1/p])$) one gets
examples of finitely presented hyperbolic groups with finitely
presented maximal subgroups of infinite index.

A similar question was suggested by Margulis and So\u{\i}fer for
the groups $\PSL_3(\Z)$. Every normal subgroup of $\PSL_3(\Z)$ is
of finite index and actually contains a principle congruence
subgroup, hence the normal topology coincides with the profinite
topology and with the congruence topology. Thus being prodense
simply means being profinitely dense, and maximal subgroups of
infinite index are automatically prodense.

\begin{question} (Margulis-So\u{\i}fer \cite{MS:first})
Does $\SL_3(\Z)$ contain a finitely generated maximal subgroup of infinite index?
\end{question}

Finally, we suppose that the methods introduced in the current paper are applicable in many other
settings.
\begin{problem} Find other settings to investigate primitivity using the machinery developed
in this paper.
\end{problem}
A natural example is the group $\Out(F_n)$, of outer automorphisms of the free group and its
subgroup. Even more interesting challenge would be to say something about countable groups acting
discretely on locally compact $\CAT(0)$ spaces. This includes linear groups as a special case, and
therefore requires in particular some geometric replacement for the condition we gave in terms of
Zariski closure.

Finally it is interesting to find geometric settings where there are few or no maximal subgroups
of infinite index. Results in this direction can be found in a recent paper by Pervova,
\cite{Perova:Grigorchuk}. In the setting of branch groups, she gives examples of groups with no
maximal subgroups of infinite index, including many Grigorchuk groups as well as groups
constructed by Gupta and Sidki.

%======================================================================================
\appendix
\section{Quasiprimitive actions in the presence of elementwise commuting normal subgroups}
\label{app:banal} We give a short analysis of primitive actions of banal groups. The methods are
similar to those appearing in the proof of O'Nan-Scott theorem (see \cite{Praeger:quasiprimitve}).
We include a short proof, because we did not find the precise statements that we need in the
literature. The results we obtain are coarser than those obtained in a classical O'Nan-Scott type
theorem. This is because we cannot assume that the commuting normal subgroups $M,N$ contain
minimal normal subgroups. The following analysis enables us to complete the proof of Theorem
\ref{thm:Linear}.

\begin{proposition} \label{prop:MN}
Let $\Gamma$ be a banal quasiprimitive group. In other words $\Gamma$ is quasiprimitive, and there
exist two non-trivial normal subgroups $\trivgp \ne N,M \lhd \Gamma$ such that  $[N,M] = \trivgp$.
Then the following hold:
\begin{enumerate}
\item \label{itm:min_N} $M$ and $N$ are the unique minimal normal subgroups of $\Gamma$. In
   particular either $M = N$ or $M \cap N = \trivgp$.
\item \label{itm:unique} The faithful quasiprimitive action of
$\Gamma$ is unique up to isomorphism
   of actions. This action is primitive.
\item \label{itm:pfaffine} If $M = N$, then $\Gamma$ is primitive
of affine type as.
   Furthermore if $\Gamma$ is finitely\footnote{recall that we exclude the case of finite groups} generated then $M$ cannot be finite dimensional as a vector space.
\item \label{itm:pfdiagonal} If $M \cap N = \trivgp$ then $\Gamma$
is primitive of diagonal type.
\end{enumerate}
\end{proposition}

\begin{proof}
Let $\Gamma \circlearrowleft \Omega$ be a faithful quasiprimitive action, $\omega \in \Omega$ and
$\Delta = \Gamma_{\omega}$. Both $M,N$ act transitively. Since $M_w$ commutes elementwise with the
transitive group $N$ it must fix $\Omega$ pointwise, hence $M$ is regular (i.e. acts freely and
transitively) and $\Gamma = \Delta \ltimes M$. Any non-trivial subgroup $M_1 < M$ which is normal
in $\Gamma$ will still be regular by the same argument, therefore $M$ is a minimal normal
subgroup. Assume that $N\neq M$. If $L$ is a minimal normal subgroup different from $M$ then $L
\cap M = \trivgp$, but then $[LN,M] = \trivgp$ so $LN$ is regular and hence $L = N$. This proves
(\ref{itm:min_N}).

Any non-trivial factor action will still be quasiprimitive and
faithful, therefore $N$ will act regularly on any such factor.
This implies that there are no non-trivial factors, so the given
action is primitive. Furthermore this primitive action is uniquely
determined. In fact since $M$ is regular we can identify $\Omega$
with $M$ using the orbit map $m \mapsto m \cdot \omega$. The
action on $\Omega$ can then be identified with the affine action
of $\Gamma$ on $M$. This proves (\ref{itm:unique}). Since the
action is primitive there are no $\Delta$ invariant subgroups of
$M$. In particular $M$ has to be characteristically simple. Note
that even for linear groups, this does not imply that $M$ is a
product of simple groups (see for example \cite{Wilson:CS}).

Assume that $M = N$. Then $M$ is abelian. Since $M$ is characteristically simple, for every prime
$p$ it is either divisible by $p$ or of exponent $p$. If $M$ has exponent $p$ for some prime $p$
then it is a vector space over $\F_p$. If $M$ does not have exponent $p$ for any $p$, then it must
be a divisible torsion free Abelian group, i.e. a vector space over $\Q$. In any case $M$ can be
identified as the additive group of a vector space over a prime field $F$. The action of $\Delta$
gives rise to a linear representation $\Delta \arrow \GL(M)$ over $F$. Finally we claim that if
$\Gamma$ is finitely generated then $M$ cannot be finite dimensional as a vector space. This is
clear if $\F = \F_p$. If $\F = \Q$ then $\Delta$ actually maps into $\GL_n(R)$ where $R$ is some
finitely generated subring of $\Q$. But then $R^n$ would be a $\Delta$ invariant subgroup of $M$
in contradiction to primitivity. This proves (3).

Now assume that $M \cap N = \trivgp$. Since both $M$ and $N$ are regular, we can define a
bijection between the two groups $\phi:M \arrow N$ by the requirement
$$
 \phi(m) \cdot \omega = m^{-1} \cdot \omega.
$$
We can easily verify that $\phi$ is a group homomorphism. Indeed assume that $\phi(m) = n$ and
$\phi(m') = n'$, using the fact that the two groups commute we can write:
\begin{eqnarray*}
(mm')^{-1} \cdot \omega & = & (m')^{-1} \cdot (m^{-1} \cdot \omega)  = (m')^{-1} \cdot (n \cdot \omega) \\
                       & = & ((m')^{-1} n) \cdot \omega           = (n (m')^{-1}) \cdot \omega \\
                       & = & n \cdot ((m')^{-1} \cdot \omega)     = n \cdot (n' \cdot \omega)\\
                       & = & (n n') \cdot \omega
\end{eqnarray*}
which implies that $\phi(mm') = nn'$. It is even easier to verify that $\phi(m^{-1}) =
\phi(m)^{-1}$. Thus $\phi$ is a group isomorphism. Denote $\Delta_0 = \{(m,\phi(m)): m \in
M\}=\Delta \cap (M \times N)$. The actions of $\Delta$ on $\Delta_0$ by conjugation is isomorphic
to its action on $M$ by conjugation, which in turn is isomorphic to the given action on $\Omega$.
In particular this action is faithful and we can identify $\Delta_0 \le \Delta \le
\Aut(\Delta_0)$. Under this identification $\Delta_0$ coincides with the group of inner
automorphisms. This establishes (4) and completes the proof of the proposition.
\end{proof}

\begin{proof}[Proof of Theorem \ref{thm:Linear}]
If $\Gamma$ is not banal then by Section \ref{sec:Linear} it is quasiprimitive if and only if it
satisfies the linear conditions for quasiprimitivity. If $\gC$ is banal then it is of affine or
diagonal type as follows from Proposition \ref{prop:MN}. It is also easy to see that groups of
affine or diagonal type are indeed primitive and that the three conditions are mutually exclusive:
if $\Gamma$ satisfies the linear conditions for primitivity then the Zariski closure of every
normal subgroup of $\Gamma$ contains the identity component of the Zariski closure of $\Gamma$
which is semisimple. The existence of elementwise commuting normal subgroups, will then imply that
the identity component must be commutative, which is absurd.
\end{proof}

We shall conclude this appendix with two basic results.

\begin{lemma} \label{lem:finite_normal}
An infinite quasi-primitive group does not contain a finite non-trivial normal subgroup.
\end{lemma}
\begin{proof}
Let $\Gamma \circlearrowleft X$ be a faithful quasi-primitive action of an infinite group. Then by
faithfulness $X$ must be infinite, and as any non-trivial normal subgroup $N \lhd \Gamma$ acts
transitively on $X$ we get $|N|\geq |X|$.
\end{proof}

\begin{corollary} \label{cor:RF_non_banal}
Assume that $\Gamma$ is a residually finite quasiprimitive group, then $\Gamma$ is never banal. In
particular the normal topology is well defined on $\Gamma$ and $\Gamma$ contains no non-trivial
normal solvable subgroups.
\end{corollary}

\begin{proof}
By Proposition \ref{prop:MN} a banal group always has a minimal
normal subgroup, which is necessarily infinite by Lemma
\ref{lem:finite_normal}. But this is impossible in a residually
finite group. If $\Gamma$ has a normal solvable subgroup $S$, then
it also has a normal abelian subgroup $A \lhd \Gamma$ -- the
smallest non-trivial group in the derived series of the solvable
group and therefore $\Gamma$ is banal.
\end{proof}

%============================================================================================

% BibTeX users please use
\bibliographystyle{alpha}
\bibliography{../tex_utils/yair}
\end{document}